\documentclass[preprint,12pt]{elsarticle}
\usepackage{lineno}
\usepackage[breaklinks]{hyperref}
\usepackage{amssymb}
\usepackage{subcaption}
\usepackage{ntheorem}
\usepackage{amsmath}
\usepackage{multirow,booktabs}
\theorembodyfont{\upshape}
\theoremseparator{.}
\newproof{pot4}{\bf Proof of Theorem \ref{theorem4}}
\newproof{pot5}{\bf Proof of Theorem \ref{theorem5}}
\newproof{pot6}{\bf Proof of Theorem \ref{theorem6}}
\newtheorem{thm}{Theorem}[section]
\newdefinition{remark}{Remark}[section]
\newtheorem{proposition}{Proposition}[section]
\newtheorem{lemma}{Lemma}[section]
\newtheorem{definition}{Definition}[section]

\newtheorem*{pro}{Proof}
\numberwithin{equation}{section}

\begin{document}
\begin{frontmatter}
\title{Global boundedness and Allee effect for a nonlocal time fractional p-Laplacian reaction-diffusion  equation \tnoteref{t1}}
\tnotetext[t1]{This work is supported by the State Key Program of National Natural Science of China under Grant No.91324201. This work is also supported by the Fundamental Research Funds for the Central Universities of China under Grant 2018IB017, Equipment Pre-Research Ministry of Education Joint Fund Grant 6141A02033703 and the Natural Science Foundation of Hubei Province of China under Grant 2014CFB865.}

\author[mymainaddress]{Hui Zhan}
\ead{2432593867@qq.com}
\author[mymainaddress]{Fei Gao\corref{mycorrespondingauthor}}
\ead{gaof@whut.edu.cn}
\author[mymainaddress]{Liujie Guo}
\ead{2252987468@qq.com}

\cortext[mycorrespondingauthor]{Corresponding author.}

\address[mymainaddress]{Department of Mathematics and Center for Mathematical Sciences, Wuhan University of Technology, Wuhan, 430070, China}

\begin{abstract}
The global boundedness and asymptotic behavior are investigated for the solutions of a nonlocal time fractional p-Laplacian reaction-diffusion  equation (NTFPLRDE)
 $$ \frac{\partial^{\alpha }u}{\partial t^{\alpha }}=\Delta_{p} u+\mu u^{2}(1-kJ*u) -\gamma u, \qquad(x,t)\in\mathbb{R}^{N}\times(0,+\infty)$$
	with $0<\alpha <1,\beta, \mu ,k>0,N\leq 2$ and $\Delta_{p}u =div(\left| \bigtriangledown  u \right|^{p-2}\bigtriangledown u)$. Under appropriate assumptions on $J$ and the conditions of $1<p<2$, it is proved that for any nonnegative and bounded initial conditions, the problem has a global bounded classical solution if $k^{*}=0$ for $N=1$ or $k^{*}=(\mu C^{2}_{GN}+1)\eta^{-1}$ for $N=2$, where $C_{GN}$ is the constant in Gagliardo-Nirenberg inequality. With further assumptions on the initial datum, for small $\mu$ values, the solution is shown to converge to $0$ exponentially or locally uniformly as $t \rightarrow \infty$, which is referred as the Allee effect in sense of Caputo derivative. Moreover, under the condition of $J \equiv 1$, it is proved that the nonlinear NTFPLRDE has a global bounded solution in any dimensional space with the nonlinear p-Laplacian diffusion terms $\Delta_{p} u^{m}\, (2-\frac{2}{N}< m\leq 3)$.
\end{abstract}

\end{frontmatter}

\section{Introduction}
 In this work we study the nonlocal  time fractional p-Laplacian reaction-diffusion equation
\begin{eqnarray}
\frac{\partial^{\alpha }u}{\partial t^{\alpha }}&=&\Delta_{p} u+\mu u^{2}(1-kJ*u) -\gamma u,\quad(x,t)\in\mathbb{R}^{N}\times(0,+\infty),\label{1}\\
u(x,0)&=&u_{0}(x),\quad x\in \mathbb{R}^{N}, \label{2}
	\end{eqnarray}
	with $0< \alpha <1,N\leq 2,\mu ,k,\gamma >0$ and $\Delta_{p}u =div(\left| \bigtriangledown  u \right|^{p-2}\bigtriangledown u),1<p<2$. By \cite{0Global} , Here $J(x)$ is a competition kernel  with
\begin{equation}\label{3}
  0\leq J \in L^{1}(\mathbb{R}^{N}),\int_{\mathbb{R}^{N}}J(x)dx=1,\inf\limits_{B(0,\delta_{0})} J> \eta
  \end{equation}
   for some $\delta _{0}>0,\eta>0$. where $B(0,\delta_{0})=(-\delta_{0},\delta_{0})^{N}$, $$J*u(x,t)=\int_{\mathbb{R}^{N}}J(x-y)u(y,t)dy.$$
And $\partial _{t}^{\alpha }$ denote the left Caputo fractional derivative that is usually defined by the formula
\begin{equation}\label{4}
  \partial _{t}^{\alpha }u(x,t)=\frac{1}{\Gamma (1-\alpha )}\int_{0}^{t}\frac{\partial u}{\partial s}(x,s)(t-s)^{-\alpha }ds,\quad 0<\alpha<1\quad
\end{equation}
in the formula \ref{4}, the left Caputo fractional derivative $\partial _{t}^{\alpha }u$ is a derivative of the order $\alpha, 0<\alpha<1$ by \cite{2018M}. Here, $F(u,J*u)=\mu u^{2}(1-kJ*u)-\gamma u$ is represented as sexual reproduction in the population dynamics system, where $-\gamma u$ is population mortality. In particular, Eq. \ref{1} is a possible model for the diffusion system of some biological species with human-controlled distribution where $u(x,t)$ represents the density of the species at position $x$ and time $t$, $div(\left| \bigtriangledown  u \right|^{p-2}\bigtriangledown u)$ portrays the mutation (which we view as a spreading of the characteristic), $J*u$ describes the resource consumption of individuals in a certain area around point $x$ at that point. This type of nonlocal term is used in biological phenomenon models such as emergence and evolution of biological species and the process of speciation in\cite{0Global}. Thus, we can find different types of nonlocal items in the nonlocal reaction-diffusion equation have different meaning and effects.

From our previous extensive research on the nonlocal reaction-diffusion equation, we found that nonlocality has a profound impact on the reaction-diffusion equation. Compared with the traditional reaction-diffusion equation, the nonlocal reaction-diffusion equation has new mathematical charateristics and richer nonlinear dynamics. Therefore, many researches are devoted to the nonlocal reaction-diffusion equation. In \cite{2015Preface}, the proposed nonlocal reaction-diffusion equation describes the population intensive distribution under the assumption of non-local sources. Its nonlocalization form is $F(u,J(u))=au^{k}(1-J(u))-\sigma u,J(u)=\int_{-\infty }^{+\infty}\phi (x-y)u(y,t)dy$. Here $k$ is a positive integer, $k=1$ corresponds to asexual and $k=2$ to sexual reproduction. we can find that the nonlocal items in the nonlinear nonlocal reaction diffusion equation also have corresponding meanings and effects. Through the investigation of \cite{2017Wavefronts}, The nonlocal term $u^{2}(1-J_{\sigma }*u)-du$ consists of the reproduction which is proportional to square of the density, the available resources and the mortality, where $-du$ is the mortality term and $d$ is the death rate. In\cite{2019Global}, its nonlocal form$\mu u^{\alpha }(1-kJ*u^{\beta }),\alpha ,\beta \geq 1$ is used describe the selection process in the mutation process, where $J*u^{\beta }(x,t)=\int _{\mathbb{R}^{N}}J(x-y)u^{\beta }(y,t)dy$ is seen to characterize the evolution of a population of density $u$. However, through many studies have shown that not only the nonlocal item has different meaning and effects, but also different nonlocal reaction-diffusion equation have different roles in various fields.

In the past few years, the nonlocal reaction-diffusion equation has been widely used in biological models, medicine and other field. For example, In\cite{2019Mathematical}, we can know characterize a cell’s motility bias according to its interactions with other cellular and acellular components in its vicinity (e.g. cell–cell and cell–tissue adhesions, nonlocal chemotaxis). This type of equation is not only used as a study of cell model research, but also can be used to explain uninfected viruses and cancer cells by \cite{2020Global}. From \cite{2011Population}, In order to predict biological populations, local competition and mutation, It proposed to solve the non-local reaction diffusion equation with numerical experiments.In addition to the above applications, the nonlocal reaction-diffusion equation will also be used by modeling emergence and evolution of a biology species in \cite{2009The}\cite{2014Elliptic}. Besides, in\cite{0Global}, it denotes by $u(x,t$the density of individuals having phenotype $x$ at time $t$ and formulate the dynamics of the population density as the following nonlinear nonlocal reaction–diffusion equation
\begin{eqnarray*}
\frac{\partial u}{\partial t}&=&\Delta u+\mu u^{2}(1-kJ*u) -\gamma u,\quad(x,t)\in\mathbb{R}^{N}\times(0,+\infty),\\
u(x,0)&=&u_{0}(x),\quad x\in \mathbb{R}^{N},
	\end{eqnarray*}
where $J*u(x,t)=\int_{\mathbb{R}^{N}}J(x-y)u(y,t)dy$. From the above studies, we studied the presence, stability and numerical experiment calculation of nonlocal reaction-diffusion equations for different applications. In order to better understand the changes in cancer cells in cell population in cancer diseases, the researchers put forward the time fractional reaction-diffusion equation.

The study of cancer cells has made a corresponding study for a long time. The earliest and most prosperous methods have been mentioned in \cite{2018M}, which mainly uses realistic data and cell model to introduce the observed phenomenon.In addition, this type of model study has significant contributions in the development of anticancer therapy. A large number of Mathematicalical tools have introduced different types of model biology in cancer. It can be classified as normal differential equations, partial differential equations, random processes, honeycomb automators and agents; see \cite{T2003A}. In \cite{2021zhan}, The following time fractional reaction-diffusion equation is used to model the biological phenomena such as emergence and evolution of biological species and the process of speciation
\begin{eqnarray*}
\frac{\partial^{\alpha }u}{\partial t^{\alpha }}&=&\Delta u+\mu u^{2}(1-kJ*u) -\gamma u,\quad(x,t)\in\mathbb{R}^{N}\times(0,+\infty),\\
u(x,0)&=&u_{0}(x),\quad x\in \mathbb{R}^{N}.
	\end{eqnarray*}
Where $J*u$ describes the resource consumption of individuals in a certain area around point $x$ at that point. And definition of $\partial _{t}^{\alpha }u$ is the same as \ref{4}. Through the study of \cite{2017Global}, it proposed the definition of Caputo derivatives in the fractional Sobolev space, and from the anoretical investigation from the operator. Specific expression forms of Caputo derivatives are shown in \ref{4}. In \cite{2015On}, it proved the existence of global bounded solutions by relying on a maximal regularity result for fractional linear reaction–diffusion equations that has been derived by Bajlekova.Also, the Caputo derivative of this equation is the same as the definition of \cite{2017Global}. Not only that, but also the same Caputo derivative in \cite{Kilbas2006}\cite{2014Basic}. Therefore, this paper will use the above-mentioned Caputo derivative to resolve the global boundaries and Allee effect of the time fractional reaction diffusion p-Laplacian equation.

In recent years, In order to better study the moderate mode of demographic dynamics and biological science, some researchers have proposed p-Laplacian reaction-diffusion equation. Thought \cite{2021Dynamic}, in order to  prove that the solution blows up at finite time, for some initial data and additional energy type conditions, it  studys the following p-Laplacian reaction–diffusion equation
\begin{equation*}
    u_{t}=div(\left | \bigtriangledown u \right |^{p-2}\bigtriangledown u)+k(t)f(u)
\end{equation*}
where $p\geq 2$ is real number, $\Omega$ is a bounded domain in $\mathbb{R}^{N}(n\geq2)$, $f$ is the source function and the time switching function $k(t)$ is nonnegative differential.The blow-up phenomena of solutions to various nonlinear problems, particularly for hyperbolic and parabolic systems, have received considerable attention in the recent literature. For example, In \cite{2014qqq}, it studied the following nonlinear p-Laplacian reaction–diffusion equation with time-independent source function
\begin{equation*}
    u_{t}=div(\left | \bigtriangledown u \right |^{p-2}\bigtriangledown u)-f(u),\qquad x\in \Omega ,t\in (0,t^{*})
\end{equation*}
subject to nonlinear boundary conduction function and initial data. In \cite{2010Periodic}, in this population dynamics model, the species is restricted to the bounded heterogeneous environment $\Omega$  whose boundary is prohibitive to the species. Thus, it leads us to the following periodic reaction -diffusion problem with p-Laplacian
\begin{align*}
    &u_{t}-div(\left | \bigtriangledown u \right |^{p-2}\bigtriangledown u)=u^{\alpha }(m(x,t)-b(x,t)u^{\beta }),\quad (x,t)\in \Omega \times \mathbb{R}^{N},\\
&u(x,t)=0,\quad (x,t)\in \partial \Omega \times \mathbb{R}^{N},\\
&u(x,0)=u_{0}(x),\quad x\in\Omega ,
\end{align*}
where $p>2,1\leq \alpha<p-1,\beta>0$, the functions $m$ and $b$ are continuous and T-periodic, and 
$\Omega$ is a bounded domain of $\mathbb{R}^{N}$ with smooth boundary $\partial \Omega $. In summary, a variety of problems in mathematical models can be solved using the P-Laplacian reaction–diffusion equation in various fields, such as establishing model in population
dynamics and biological sciences, the blow-up phenomena and population dynamics, so on.

In addition to the equation of the p-Laplace type is generated in many applications in physical and biological fields, the solutions of this type equation also have different properties, such as existence,stability, blow-up, etc. In recent years, the existence of the p-Laplacian equation has been widely concerned \cite{2004Theaa}\cite{2016Existence}, which is mainly carried out under different boundary conditions. In \cite{2013Stability}, the  mathematical model mainly consider the stability of blowup of solutions for the p-Laplace equation with nonlinear source. By \cite{2003Blow}, the following blow-up criteria for weak solutions of the Dirichlet problem
$$u_{t}=\bigtriangledown (\left | \bigtriangledown u \right |^{p-2}\bigtriangledown u)+\lambda \left | u \right |^{q-2}u,\quad in\quad \Omega _{T}$$
where $p>1$, which shows the blow-up properties of the above semilinnear heat equation. Many authors have discussed the blow-up phenomena of p-Laplacian parabolic problems. In \cite{Ding2018Global}, it studys the global existence and blow-up results for the p-Laplacian parabolic problems by relying mainly on constructing some auxiliary functions and using the parabolic maximum principles and the differential inequality technique. Therefore, the solitation of p-Laplacian equation is different, and the above properties are different from the above properties. In addition to the existence, stability,  and blow-upof the p-Laplacian equation,This type of p-Laplacian equation also has extinctive and non-extinction.

In the last decades, many researchers devoted to the study of extinction and non-extinction of solutions for nonlinear parabolic equations with nonlocal terms. In \cite{2003Global}, the following equation studied the homogeneous Dirichlet boundary value problem for the degenerate parabolic equation with nonlocal source term
$$u_{t}-div(\left | \bigtriangledown u \right |^{p-2}\bigtriangledown u)=\int _{\Omega }u^{q}(x,t)dx,\quad x\in\Omega ,t>0$$
where $\Omega$ is a bounded domain in $\mathbb{R}^{N}(N\geq1),p>2,q\geq 1$.The extinction phenomenon of solutions for these equations also makes some progress. thought the study of \cite{2010Extinction}\cite{2001THE}\cite{2001Extinctionq}, They obtained sufficient conditions about the extinction and non-extinction of solutions for the p-Laplacian equation by the upper and lower solutions methods. In references \cite{2012Extinction}, the following equation study the extinction, non-extinction and decay estimates of non-negative nontrivial weak solutions of the initial-boundary value problem for the p-Laplacian equation with nonlocal nonlinear source and interior linear absorption.
$$u_{t}=div(\left | \bigtriangledown u \right |^{p-2}\bigtriangledown u)-\lambda\int _{\Omega }u^{q}(x,t)dx-ku,\quad x\in\Omega ,t>0$$
where $1<p<2,k,q,\lambda>0,\Omega \subset \mathbb{R}^{N} (N\geq 1)$ is a bounded domain with smooth boundary. This article will use this definition based on the definition of P-Laplace in the above equation. In the following we will switch our view point to investigate the global boundedness and Allee effect of the solutions for the corresponding nonlocal time fractional reaction-diffusion p-
Laplacian equation.

In the paper, we put forward the concept of the nonlocal time fractional reaction-diffusion p-
Laplacian equation, in which sufficient conditions are given for the global boundedness of its solution, asymptoticity and global boundedness under nonlinear conditions. First, we mainly take advantage of the inequality of Bernoulli, Sobolev embedding inequality and fractional Duhamel's formula \cite{2018Duhamel} to prove that theorem \ref{theorem4}. Secondly, in order to prove the Allee effect of the nonlocal time fractional reaction-diffusion p-Laplacian equation, our main critical use the fractional differential equations. Fractional differential equations which are differential equations of arbitrary order have attracted considerable attention recently \cite{2009Existence}\cite{2002Analytical}.
However, we find that the existence of Laplace transform is taken for granted in some papers to solve fractional differential equations (see \cite{2007Global}\cite{2008Stability}).In \cite{2011Laplace}, Giving a sufficient condition to guarantee the rationality of solving constant coefficient fractional differential equations by the Laplace transform method. Finally, we will mainly use Gagliardo-Nirenberg inequality to prove that nonlinear nonlocal time fractional reaction-diffusion p-
Laplacian equation in $L^{r}$ estimates and $L^{\infty}$ estimates, respectively, with global boundary.
Next, we study the global boundedness and asymptotic for the solution of nonlocal time fractional reaction-diffusion p-Laplacian equation of \eqref{1}-\eqref{4}. The main results of the paper are the following.

\begin{thm}\label{theorem4}
  Suppose \eqref{3}-\eqref{4} holds, $0\leq u_{0}\in L^{\infty }(\mathbb{R}^{N})$ ,$1<p<2$ and $T$ is constant. Denote $k^{*}=0$ for $N=1$ and $k^{*}=(\mu C_{GN}^{2}+1)\eta ^{-1}$ for $N=2$, $C_{GN}$ is the constant appears in Gagliardo-Nirenberg inequality in Lemma \ref{lemma2.1}, then for any $k>k^{*}$, the nonnegative solution of \eqref{1}-\eqref{2} exists and is globally bounded in time, that is, there exist
  \begin{equation}
 \nonumber   K=\begin{cases}
 K(\Arrowvert u_{0}\Arrowvert_{L^{\infty}(\mathbb{R}^{N})},\mu,p,\eta,k,C_{GN},T^{\alpha}),&  N=1,\\
 K(\Arrowvert u_{0}\Arrowvert_{L^{\infty}(\mathbb{R}^{N})},\mu,p,C_{GN},T^{\alpha }),&  N=2,
      \end{cases}
  \end{equation}
  such that
\begin{equation}\label{7}
     0\leq u(x,t)\leq K,\qquad\forall (x,t)\in \mathbb{R} ^{N}\times [0,\infty).
  \end{equation}
\end{thm}

\begin{thm}\label{theorem5}
  Denote $u(x,t)$ the globally bounded solution of \eqref{1}-\eqref{2}.
\begin{enumerate}[\bf (i).]
\item For any $\gamma>0$, there exist $\mu ^{*}>0$ and $m^{*}>0$ such that for\\ $\mu \in (0,\mu ^{*})$ and $\left \| u_{0} \right \|_{L^{\infty }(\mathbb{R}^{N})}< m^{*}$, we have $$\left \| u(x,t) \right \|_{L^{\infty }(\mathbb{R}^{N}\times [0,+\infty ))}< \frac{\tau }{\mu },$$ and thus $$\left \| u(x,t) \right \|_{L^{\infty }(\mathbb{R}^{N})}\leq \left \| u_{0} \right \|_{L^{\infty }(\mathbb{R}^{N})}e^{\left ( -\sigma  \right )^{\frac{1}{\alpha }}t},$$ for all $t>0,0<\alpha<1$ with $\sigma :=\gamma -\mu \left \| u(x,t) \right \|_{L^{\infty }(\mathbb{R}^{N}\times [0,+\infty ))}>0.$

\item If $0< \gamma < \frac{\mu }{4k} $, there exist $\mu ^{**}>0$ and $m^{**}>0$ such that for $\mu \in (0,\mu ^{**})$ and $\left \| u_{0} \right \|_{L^{\infty }(\mathbb{R}^{N})}< m^{**}$, we have $$\left \| u(x,t) \right \|_{L^{\infty }(\mathbb{R}^{N}\times [0,+\infty ))}< a,$$ and thus $$\lim_{t\rightarrow \infty }u(x,t)=0$$ locally uniformly in $\mathbb{R}^{N}$.
  \end{enumerate}
\end{thm}

\begin{thm}\label{theorem6}
 \begin{eqnarray}
  \frac{\partial^{\alpha } u}{\partial t^{\alpha }}&=&\Delta_{p} u^{m}+u^{2}(1-\int _{\mathbb{R}^{N}}udx)-u(x,t)  \label{1.1.4}\\
u(x,0)&=&u_{0}(x) \qquad x \in \mathbb{R}^{N}\label{1.1.5}
 \end{eqnarray}
with $ (x,t)\in \mathbb{R}^{N}\times(0,+\infty),1<p<2$. If $2-\frac{2}{N}< m\leq 3,0< \alpha <1$, then the solution to any initial value $0\leq u_{0}\in L^{\infty }(\mathbb{R}^{N})$ exists and \eqref{1.1.4}-\eqref{1.1.5} is globally bounded, exist $M>0$ such that$$0\leq u(x,t)\leq M, \qquad (x,t)\in \mathbb{R}^{N}\times [0,+\infty).$$
	\end{thm}
	
\section{Global bounded of solutions for a NTFPLRDE}

\begin{lemma}\textsuperscript{\cite{Alikhanov2010A}}\label{lemma23}
The left Caputo derivative with respect to time $t$ of $v$ is defined by \eqref{4} and absolutely continuous on $[0,T]$, one has the inequality
  \begin{equation}\label{888}
    v(t)\partial _{t}^{\alpha }v(t)\geq \frac{1}{2}\partial _{t}^{\alpha }v^{2}(t),\quad 0<\alpha <1.
  \end{equation}
\end{lemma}

\begin{lemma}\textsuperscript{\cite{2021zhan}}\label{lemma2}
Suppose $u:[0,\infty )\times\mathbb{R}^{n}\rightarrow \mathbb{R}$,the left Caputo fractional derivative with respect to time $t$ of $u$ is defined by \eqref{4}. Then there is
\begin{equation}\label{6}
  \int _{B(x,\delta )}\partial _{t}^{\alpha }udy=\partial _{t}^{\alpha }\int _{B(x,\delta )}udy.
\end{equation}
	\end{lemma}
	
\begin{definition}\textsuperscript{\cite{Kilbas2006}}
Assume that $X$ is a Banach space and let $u:\left [ 0,T \right ]\rightarrow X$. The Caputo fractional derivative operators of $u$ for order $\alpha \in \mathbb{C},(Re(\alpha )>0)$ are defined by
\begin{eqnarray}
_{0}^{C}\textrm{D}_{t}^{\alpha }u(t)&=&\frac{1}{\Gamma (1-\alpha )}\int_{0}^{t}(t-s)^{-\alpha }\frac{d}{ds}u(s)ds,\label{32}\\
\nonumber_{t}^{C}\textrm{D}_{T}^{\alpha }u(t)&=&\frac{-1}{\Gamma (1-\alpha )}\int_{t}^{T}(s-T)^{-\alpha }\frac{d}{ds}u(s)ds,
\end{eqnarray}
 where $\Gamma (1-\alpha )$ is the Gamma function. The above integrals are called the left-sided and the right-sided the Caputo fractional derivatives
\end{definition}

\begin{lemma}\textsuperscript{\cite{2021zhan}}\label{lemma1}
Suppose $u:[0,\infty )\times \mathbb{R}^{n}\rightarrow \mathbb{R}$ and $n\geq 2$,the Caputo fractional derivative with respect to time $t$ of $u$ is defined by \eqref{32}. Then there is
\begin{equation*}
   u^{n-1}(_{0}^{C}\textrm{D}_{t}^{\alpha }u) \geq \frac{1}{n}(_{0}^{C}\textrm{D}_{t}^{\alpha }u^{n}).
\end{equation*}
\end{lemma}

\begin{lemma}\textsuperscript{\cite{2021zhan}}\label{lemma3}
  Suppose $u:[0,\infty )\times\mathbb{R}^{n}\rightarrow \mathbb{R}$,the Caputo fractional derivative with respect to time $t$ of $u$ is defined by \eqref{32}. Then there is
\begin{equation*}
  \int _{B(x,\delta )}(_{0}^{C}\textrm{D}_{t}^{\alpha }u)dy=_{0}^{C}\textrm{D}_{t}^{\alpha }\int _{B(x,\delta )}udy.
\end{equation*}
\end{lemma}

\begin{lemma}\textsuperscript{\cite{2017GagliardO}}\label{lemma99}
  Assume $1\leq p<n$, then $\forall f\in C_{c}^{\infty }(\mathbb{R}^{N})$, there is 
  \begin{equation}\label{99}
      \left \| f \right \|_{L^{q}(\mathbb{R}^{N})}\leq C_{1}\left \| \bigtriangledown f \right \|_{L^{p}(\mathbb{R}^{N})}.
  \end{equation}
  Which is $q=\frac{np}{n-p}$, and $C_{1}$ only dependent on $p,n$.
\end{lemma}

\begin{lemma}\textsuperscript{\cite{2021zhan}}\label{corollary1}
  Suppose that a nonnegative function $y(t)\geq 0$  satisfies
  \begin{equation}
    \nonumber _{0}^{C}\textrm{D}_{t}^{\alpha }y(t)+c_{1}y(t)\leq b
  \end{equation}
for almost all $t\in [0,T]$, where $b,c_{1}>0$ are all contants. then
\begin{equation*}
   y(t)\leq y(0)+\frac{bT^{\alpha }}{\alpha \Gamma (\alpha )}.
  \end{equation*}
\end{lemma}

\begin{lemma}
  Let $y(t)$ be a non-negative absolutely continuous function on $[0,+\infty)$ satisfying
  $$_{0}^{C}\textrm{D}_{\alpha }^{t}y(t)+C_{1}y^{k}(t)+C_{2}y(t)\leq C,\quad t\in [0,T]$$
where $C_{1},C_{2}>0$ are constantss and $k\in (0,1)$, then we have decay estimate
\begin{equation}
    y(t)\leq [y^{1-k}(0)+\frac{(C+C_{1}(k-1))T^{\alpha }}{\alpha \Gamma(\alpha )}]^{\frac{1}{1-k}}
\end{equation}
\end{lemma}

\begin{pro}
    Make $z(t)=y^{1-k}(t)$, then the above formula can be written
    $$_{0}^{C}\textrm{D}_{\alpha }^{t}z(t)+(1-k)C_{2}z(t)\leq (k-1)C_{1}+C.$$
By Lemma \ref{corollary1}, then we can obtain the following inequality
$$z(t)\leq z(0)+\frac{(C+C_{1}(k-1))T^{\alpha }}{\alpha\Gamma(\alpha )}.$$
So, We get conclusions
$$y(t)\leq [y^{1-k}(0)+\frac{(C+C_{1}(k-1))T^{\alpha }}{\alpha\Gamma(\alpha )}]^{\frac{1}{1-k}}.$$
\end{pro}

\begin{remark}
Through the study of such inequality in reference\cite{2012Extinction}\cite{2010Extinction}\cite{2001THE}\cite{1998The}, the variation of Bernoulli is made in this paper. And the proven of this inequality is also verified by using the Bernoulli equation. This plays a very important role in the latter theorem prove.
\end{remark}

Before the certification theorem \ref{theorem4}, we firstly will use fractional Doumal's  formula to get the solution of equation \eqref{1}-\eqref{2}. We define the operator $A$ on $L^{2}(\Omega )$ by
\begin{equation*}
\begin{cases}
 D(A)=\left \{ u\in H_{0}^{1}(\Omega ):\Delta_{p} u\in L^{2}(\Omega ) \right \}, \\
 Au=\Delta_{p} u,\qquad for \, u\in D(A).
\end{cases}
\end{equation*}
From equations \eqref{1}-\eqref{2}, we consider the equations

\begin{equation}\label{2.8}
\begin{cases}
  \frac{\partial ^{\alpha }}{\partial t^{\alpha }}u(t)=Au(t)+h(t), \\
 u(0)=u_{0},
\end{cases}
\end{equation}
 where $h(t)=\mu u^{2}(1-kJ*u)-\gamma u$ and $A$ is self-adjoint on Hilbert space $X$. Then it follows from spectral theorem that there exists a measure space $(\Sigma ,\mu )$ and a Borel measurable function $a$ and a unitary map $U:L^{2}(\Sigma,\mu )\rightarrow X$ such that$$U^{-1}AU=T_{a},T_{a}\varphi (\xi )=a(\xi )\varphi (\xi ),\quad \xi \in \Sigma. $$
 \begin{lemma}\label{lem3}\textsuperscript{\cite{2021zhan}}
   If $u$ satisfies Equation \eqref{2.8}, then $u$ also satisfies
   $$
   \begin{aligned}
   u(t)&=U(E_{\alpha }(a(\xi )t^{\alpha })U^{-1}u_{0})\\
   &\qquad+\int_{0}^{t}(t-s)^{\alpha -1}U(E_{\alpha ,\alpha }((t-s)^{\alpha }a(\xi ))U^{-1}h(s))ds, \quad t\in [0,T].
    \end{aligned}
   $$
 \end{lemma}
 
 \begin{remark}
 By define two maps $$\mathcal{S} _{\alpha }(t)\phi=U(E_{\alpha }(a(\xi )t^{\alpha })U^{-1}\phi),\, \mathcal{K} _{\alpha }(t)\phi=U(E_{\alpha ,\alpha }((t)^{\alpha }a(\xi ))U^{-1}\phi),\, \phi \in X,$$then through Lemma \ref{lem3}, we can get the following fractional Duhamel's formula.
   \end{remark}
 
\begin{lemma}\label{lemma2.6}\textsuperscript{\cite{2021zhan}}
  If $u\in C([0,T],X)$ satisfies Equation \eqref{1} -\eqref{2}, then $u$ satisfies the following integral equation:
   $$u(t)=\mathcal{S}_{\alpha }(t)u_{0}+\int_{0}^{t}(t-s)^{\alpha -1}\mathcal{K} _{\alpha }(t-s)h(s)ds.$$
   The solution operators $\mathcal{S}_{\alpha }(t)$ and $\mathcal{K} _{\alpha }(t)$ are defined by the functional calculus of $A$ via the Mittag-Leffler function when evaluated at $A$.
\end{lemma}	

\begin{proposition}\textsuperscript{\cite{2018Duhamel}}\label{proposition.1}
   For each fixed $t \geq 0$,$\mathcal{S}_{\alpha }(t)$ and $\mathcal{K} _{\alpha }(t)$ are linear and bounded operators, for any $\phi \in X$,

   \begin{equation}\label{9}
     \left \|\mathcal{S}_{\alpha }(t)\phi  \right \|_{X}\leq C_{4}\left \| \phi \right \|_{X},\qquad  \left \|\mathcal{K}_{\alpha }(t)\phi  \right \|_{X}\leq C_{4}\left \| \phi \right \|_{X},
   \end{equation}
where $C_{4}$ is a constant.
 \end{proposition}

\begin{lemma}\textsuperscript{\cite{1966}}\label{lemma2.1}
  Let $\Omega$ be an open subset of $\mathbb{R}^{N}$, assume that $1\leq p,q \leq \infty$ with $(N-q)p<Nq $ and $r\in(0,p)$. Then there exists constant $C_{GN}>0$ only depending on $q, r$ and $\Omega$ such that for any $u\in W^{1.q}(\Omega)\cap L^{p}(\Omega)$
  $$\int_{\Omega}u^{p}dx\leq C_{GN}(\Arrowvert\bigtriangledown u\Arrowvert^{(\lambda^{*}p)}_{L^{q}(\Omega)}\Arrowvert u\Arrowvert_{L^{p}(\Omega)}^{(1-\lambda^{*})p}+\Arrowvert u\Arrowvert_{L^{p}(\Omega)}^{p})$$
  holds with $$\lambda^{*}=\frac{\frac{N}{r}-\frac{N}{p}}{1-\frac{N}{q}+\frac{N}{r}}\in (0,1).$$
\end{lemma}

\begin{lemma}\textsuperscript{\cite{2016Partial}}\label{Young}
  $\forall a,b\geq 0$ and $\varepsilon  >0$, for $1<p,q<\infty ,\frac{1}{p}+\frac{1}{q}=1$, then there is $$a\cdot b\leq \varepsilon\frac{a^{p}}{p} +\varepsilon ^{-\frac{q}{p}}\frac{b^{q}}{q}.$$
\end{lemma}

\begin{lemma}\textsuperscript{\cite{2019Global}}\label{lemma2.8}
  Fix $x_{0}= (x_{1}^{0},\cdots,x_{N}^{0})\in \mathbb{R}^{N}$, choose $0<\delta\leq \frac{1}{2}\delta_{0}$, and denote $B(x_{0},\delta)=(x_{0}= (x_{1}^{0},\cdots,x_{N}^{0})\in \mathbb{R}^{N}|\quad \lvert x_{i}-x^{0}_{i}\lvert\leq\delta,1\leq i\leq N)$, then for any $y\in B(x,\delta)$, there is
  \begin{equation}\label{10}
  \begin{aligned}
    \int_{B(x,\delta)}\Delta u^{2}dy&=\int_{B(0,\delta)}\Delta u^{2}(y+x,t)dy\\
    &=\Delta \int_{B(0,\delta)}u^{2}(y+x,t)dy=\Delta \int_{B(x,\delta)}u^{2}dy.
    \end{aligned}
  \end{equation}
\end{lemma}

\begin{proposition}\textsuperscript{\cite{1966On}}\label{Proposition2.1}
  Assume the initial data $0\leq u_{0}\in L^{\infty}(\mathbb{R}^{N})$. Then there are a maximal existence time $T_{max}\in(0,\infty]$ and $u\in C([0,T_{max}),L^{\infty}(R^{N}))\bigcap$\\ $C^{2,1}(\mathbb{R}^{N}\times[0,T_{max}))$ such that $u$ is the unique nonnegative classical solution of \eqref{1}-\eqref{2}. Furthermore, if $T_{max}<+\infty$, then$${\lim_{t\to T_{max}}}\Arrowvert u(\cdot,t)\Arrowvert_{L^{\infty}(\mathbb{R}^{N})}=\infty.$$
\end{proposition}

\begin{pot4}
Fix $x_{0}= (x_{1}^{0},\cdots,x_{N}^{0})\in \mathbb{R}^{N}$, choose $0<\delta\leq \frac{1}{2}\delta_{0}$, and denote $B(x_{0},\delta)=(x_{0}= (x_{1}^{0},\cdots,x_{N}^{0})\in \mathbb{R}^{N}|\quad \lvert x_{i}-x^{0}_{i}\lvert\leq\delta,1\leq i\leq N)$ with $\lvert B(x_{0},\delta)\lvert=(2\delta)^{N}$. For any $x \in \mathbb{R}^{N}$,multiply \eqref{1} by $2u\varphi_{\varepsilon}$, where $\varphi_{\varepsilon}(\cdot)\in C_{0}^{\infty}(B(x,\delta))$, and $\varphi_{\varepsilon}(\cdot)\to 1$ locally uniformly in $B(x,\delta)$ as $\varepsilon\to 0$. Integrating by parts over $B(x,\delta)$, we obtain

$$
\begin{aligned}
\int_{B(x,\delta)}2u\varphi_{\varepsilon}\frac{\partial^{\alpha } u}{\partial t^{\alpha }}dy&=\int_{B(x,\delta)}2u\varphi_{\varepsilon}\Delta_{p} udy+\int_{B(x,\delta)}2u^{3}\mu\varphi_{\varepsilon}(1-kJ*u)dy\\
&\qquad-2\gamma\int _{B(x,\delta )}u^{2}\varphi_{\varepsilon} dy,
\end{aligned}
$$
from 
\begin{align*}
    &\int_{B(x,\delta)}2u\varphi_{\varepsilon}\Delta_{p} udy=\int_{B(x,\delta)}2u\varphi_{\varepsilon}div(\left | \bigtriangledown u \right |^{p-2}\bigtriangledown u)dy\\
&=\int_{B(x,\delta)}\varphi_{\varepsilon}(\left | \bigtriangledown u \right |^{p-2}\Delta u^{2}-2\left | \bigtriangledown u \right |^{p})dy.
\end{align*}

and Lemma \ref{lemma23}  and Lemma \ref{lemma2}, then we can get
$$\int_{B(x,\delta)}2u\varphi_{\varepsilon}\frac{\partial^{\alpha } u}{\partial t^{\alpha }}dy\geq\frac{\partial^{\alpha }}{\partial t^{\alpha }}\int_{B(x,\delta)}u^{2}\varphi_{\varepsilon}dy,$$
and
$$
\begin{aligned}
   &\frac{\partial^{\alpha }}{\partial t^{\alpha }}\int_{B(x,\delta)}u^{2}\varphi_{\varepsilon}dy+2\int_{B(x,\delta)}\arrowvert\bigtriangledown u\arrowvert^{p}\varphi_{\varepsilon}dy \\
   &\leq\int_{B(x,\delta)}(\left | \bigtriangledown u \right |^{p-2}\Delta u^{2}\varphi_{\varepsilon}dy+2\mu\int_{B(x,\delta)}u^{3}(1-kJ*u)\varphi_{\varepsilon}dy\\
   &-2\gamma\int_{B(x,\delta)}u^{2}\varphi_{\varepsilon}dy.
\end{aligned}
$$
Taking $\varepsilon\to 0$, we obtain
\begin{equation}\label{11}
\begin{aligned}
&\frac{\partial^{\alpha}}{\partial t^{\alpha}}\int_{B(x,\delta)}u^{2}dy+2\int_{B(x,\delta)}\arrowvert\bigtriangledown u\arrowvert^{p}dy\\
&\leq\left | \bigtriangledown u \right |^{p-2}\int_{B(x,\delta)}\Delta u^{2}dy+2\mu\int_{B(x,\delta)}u^{3}(1-kJ*u)dy-2\gamma\int_{B(x,\delta)}u^{2}dy\\
&\leq\left | \bigtriangledown u \right |^{p-2}\int_{B(x,\delta)}\Delta u^{2}dy+2\mu\int_{B(x,\delta)}u^{3}dy-2\mu\int_{B(x,\delta)}u^{3}kJ*udy\\
&\qquad-2\gamma\int_{B(x,\delta)}u^{2}dy.
\end{aligned}
\end{equation}
Knowing from \eqref{3} and Lemma \ref{lemma2.8}, we have used the fact that $\forall y\in B(x,\delta),z\in B(y,\delta)$, then $J(z-y)\geq\eta$ and $B(x,\delta)\subset B(y,2\delta)$, and then
$$
\begin{aligned}
&J*u(y,t)=\int_{B(y,\delta)}J(y-z)u(z,t)dz\geq\eta\int_{B(x,\delta)}u(y,t)dy,\\
&\int_{B(x,\delta)}\Delta u^{2}dy=\Delta \int_{B(x,\delta)}u^{2}dy.
\end{aligned}
$$
Therefore
$$
\begin{aligned}
&\frac{\partial^{\alpha }}{\partial t^{\alpha }}\int_{B(x,\delta)}u^{2}dy+2\int_{B(x,\delta)}\arrowvert\bigtriangledown u\arrowvert^{p}dy\\
&\leq\left | \bigtriangledown u \right |^{p-2}\Delta \int_{B(x,\delta)}u^{2}dy+2\mu\int_{B(x,\delta)}u^{3}dy-2\mu\eta k\int_{B(x,\delta)}u^{3}dy\int_{B(x,\delta)}udy\\
&\qquad-2\gamma\int_{B(x,\delta)}u^{2}dy.
\end{aligned}
$$
Now we proceed to estimate the term $\int_{B(x,\delta)}u^{3}dy$. Firstly using Gagliardo-Nirenberg inequality in Lemma \ref{lemma2.1} on $\Omega=B(x,\delta)$ with $p=3,q=r=2,\lambda^{*}=\frac{N}{6}$, there exists constant $C_{GN}>0$, such that
\begin{equation}\label{22}
  \int_{B(x,\delta)}u^{3}dy\leq C_{GN}(N,\delta)(\Vert\bigtriangledown u\Vert_{L^{2}(B(x,\delta))}^{\frac{N}{2}}\Vert u\Vert_{L^{2}(B(x,\delta))}^{3-\frac{N}{2}}+\Vert u\Vert_{L^{2}(B(x,\delta))}^{3}).
\end{equation}
On the one hand, by Lemma \ref{Young}, we can make $a=\Vert\bigtriangledown u\Vert_{L^{2}(B(x,\delta))}^{\frac{N}{2}}$,\\ $b=C_{GN}(N,\delta)\Vert u\Vert_{L^{2}(B(x,\delta))}^{3-\frac{N}{2}}, \varepsilon=\frac{1}{\mu},p=\frac{1}{N},q=\frac{4}{4-N}$, then obtain
\begin{equation}\label{23}
\begin{aligned}
&C_{GN}(N,\delta)\Vert\Delta u\Vert_{L^{2}(B(x,\delta))}^{\frac{N}{2}}\Vert u\Vert_{L^{2}(B(x,\delta))}^{3-\frac{N}{2}}\\
&\leq\frac{1}{\mu}\Vert \Delta u\Vert_{L^{2}(B(x,\delta))}^{2}+\mu^{\frac{N}{4-N}}C^{\frac{4}{4-N}}_{GN}(N,\delta)\Vert u\Vert_{L^{2}(B(x,\delta))}^{\frac{2(6-N)}{4-N}}.
\end{aligned}
\end{equation}
Here by Young's inequality, let $a=C_{GN}(N,\delta ),b=\left \| u \right \|_{L^{2}(B(x,\delta ))}^{3},p=\frac{2(6-N)}{3(4-N)}$, $q=\frac{2(6-N)}{N}$, then
\begin{equation}\label{24}
  C_{GN}(N,\delta )\left \| u \right \|_{L^{2}(B(x,\delta ))}^{3}\leq \left \| u \right \|_{L^{2}(B(x,\delta ))}^{\frac{2(6-N)}{4-N}}+C^{\frac{2(6-N)}{N}}(N,\delta).
\end{equation}
By interpolation inequality, we obtain
\begin{equation}\label{25}
\begin{aligned}
\left \| u \right \|_{L^{2}(B(x,\delta ))}^{\frac{2(6-N)}{4-N}}&\leq (\left \| u \right \|_{L^{1}(B(x,\delta ))}^{\frac{1}{4}}\left \| u \right \|_{L^{3}(B(x,\delta ))}^{\frac{3}{4}})^{\frac{2(6-N)}{4-N}}\\
&=(\int _{B(x,\delta)}u^{3}dy\int _{B(x,\delta )}udy)^{\frac{6-N}{2(4-N)}}.
\end{aligned}
\end{equation}
Combing \eqref{22}-\eqref{25}, we obtain
\begin{equation}
  \begin{aligned}
  &2\mu \int _{B(x,\delta )}u^{3}dy\\
  &\leq 2\int _{B(x,\delta )}\left | \bigtriangledown u \right |^{2}dy+2\mu(1+\mu^{\frac{N}{4-N}}C_{GN}^{\frac{4}{4-N}}(N,\delta ))(\int _{B(x,\delta )}u^{3}dy\int _{B(x,\delta )}udy)^{\frac{6-N}{2(4-N)}}\\
  &\qquad+2\mu C_{GN}^{\frac{2(6-N)}{N}}(N,\delta).\label{26}
  \end{aligned}
\end{equation}
Next we consider the cases $N=1$ and $N=2$ respectively.

\textbf{Case 1.}  $N=1$. From \eqref{26}, by Lemma \ref{Young}, let $\varepsilon =\eta k,p=\frac{6}{5},q=6$, we get
\begin{equation}
  \begin{aligned}
  &2\mu \int _{B(x,\delta )}u^{3}dy\\
&\leq 2\int _{B(x,\delta )}\left | \bigtriangledown u \right |^{2}dy
+2\mu \eta k\int _{B(x,\delta )}u^{3}dy\int _{B(x,\delta )}udy\\
&\qquad+2\mu (\mu^{\frac{1}{3}}C_{GN}^{\frac{4}{3}}(1,\delta )+1)^{6}(\eta k)^{-5}+2\eta C_{GN}^{10}(1,\delta ),\label{27}
  \end{aligned}
\end{equation}
inserting \eqref{27} into \eqref{11}, we have
 \begin{equation}\label{28}
   \begin{aligned}
  &\frac{\partial^{\alpha } }{\partial t^{\alpha }}\int _{B(x,\delta )}u^{2}dy+2\int _{B(x,\delta )}\left | \bigtriangledown u \right |^{p-2}dy+2\gamma \int _{B(x,\delta )}u^{2}dy\\
&\leq\left | \bigtriangledown u \right |^{p-2} \Delta \int _{B(x,\delta )}u^{2}dy+2\int _{B(x,\delta )}\left | \bigtriangledown u \right |^{2}dy\\
&+2\mu [(\mu^{\frac{1}{3}}C_{GN}^{\frac{4}{3}}(1,\delta ))^{6}(\eta k)^{-5}+C_{GN}^{10}(1,\delta )].
  \end{aligned}
 \end{equation}
 Denote $y(t)=\int_{B(x,\delta )}u^{2}dy=\left \| u \right \|_{2}^{2}
$ and by Lemma \ref{lemma99}, We can get the following inequality
 \begin{align}\label{2.17}
     &C_{1}\left \| u \right \|_{2}^{p}\leq \left \| \bigtriangledown u \right \|_{p}^{p},\quad C_{2}\left \| u \right \|_{2}^{2}\leq \left \| \bigtriangledown u \right \|_{2}^{2},\quad p\in (1,2)
 \end{align}
 And we get the following fractional Bernui inequality
 \begin{equation*}\label{29}
   \begin{cases}
    _{0}^{C}\textrm{D}_{t}^{\alpha }y(t)+2C_{1}y^{p}(t)+(2\gamma-C_{2})y(t)\leq Q\\
y(0)=2\delta \left \| u_{0} \right \|_{L^{\infty }(R^{N})}^{2}
   \end{cases}
 \end{equation*}
where $Q_{1}=2\mu \left ( (\mu ^{\frac{1}{3}}C_{GN}^{\frac{4}{3}}(1,\delta)+1)^{6}(\eta k)^{-5}+C_{GN}^{10}(1,\delta ) \right )$, then we obtain
\begin{align}
    y(t)=\left \| u \right \|_{2}^{2}\leq [y^{1-p}(0)+\frac{(Q+(2\gamma -2C_{2}))T^{\alpha }}{\alpha \Gamma (\alpha )}]^{\frac{1}{1-p}}=M_{1}
\end{align}
 \textbf{Case 2.}$N=2$. From \eqref{26}, we have
 \begin{equation}\label{2.10}
   \begin{aligned}
   &2\mu \int _{B(x,\delta )}u^{3}dy\\
   &\leq 2\int _{B(x,\delta )}\left | \bigtriangledown u \right |^{2}dy+2\mu (1+\mu C_{GN}^{2}(2,\delta ))(\int _{B(x,\delta )}u^{3}dy\int _{B(x,\delta )}udy)\\
   &\qquad+2\mu C_{GN}^{4}(2,\delta).
   \end{aligned}
 \end{equation}
 Inserting \eqref{2.10} into \eqref{11}, for $k\geq k^{*}=\frac{\mu C_{GN}^{2}(2,\delta )+1}{\eta }$, we obtain
 \begin{align*}\label{2.11}
   &\frac{\partial^{\alpha } }{\partial t^{\alpha }}\int _{B(x,\delta )}u^{2}dy+2\gamma \int _{B(x,\delta )}u^{2}dy+2\int _{B(x,\delta )}\left | \bigtriangledown u \right |^{p}dy\\
&\leq \left | \bigtriangledown u \right |^{p}\Delta \int _{B(x,\delta )}u^{2}dy+2\mu (\mu C_{GN}^{2}(2,\delta )+1-\eta k)\int _{B(x,\delta )}u^{3}dy\int _{B(x,\delta )}udy\\
&\qquad+2\mu C_{GN}^{4}(2,\delta)+2\int _{B(x,\delta )}\left | \bigtriangledown u \right |^{2}dy\\
&\leq \left | \bigtriangledown u \right |^{p}\Delta \int _{B(x,\delta )}u^{2}dy+2\int _{B(x,\delta )}\left | \bigtriangledown u \right |^{2}dy+2\mu C_{GN}^{4}(2,\delta).
\end{align*}
 Denote $y(t)=\int_{B(x,\delta )}u^{2}dy$ and inequality \ref{2.17}, then  getting the following Bernui inequality
 \begin{equation*}\label{2.12}
   \begin{cases}
    _{0}^{C}\textrm{D}_{t}^{\alpha }y(t)+2C_{1}y^{p}(t)+(2\gamma-C_{2})y(t)\leq 2\mu C_{GN}^{4}(2,\delta),\\
   y(0)=(2\delta)^{2}\left \| u_{0} \right \|_{L^{\infty}(\mathbb{R}^{N})}^{2}.
   \end{cases}
 \end{equation*}
  Then, we obtain 
  \begin{align}
    y(t)=\left \| u \right \|_{2}^{2}\leq [y^{1-p}(0)+\frac{(2\mu C_{GN}^{4}(2,\delta)+(2\gamma -2C_{2}))T^{\alpha }}{\alpha \Gamma (\alpha )}]^{\frac{1}{1-p}}=M_{2}
\end{align}
 In conclusion, for any $(x,t)\in  \mathbb{R}^{N}\times [0,T_{max})$, we have
 \begin{equation}\label{2.13}
   \left \| u \right \|_{L^{2}(B(x,\delta ))}\leq M=\begin{cases}
	\sqrt{M_{1}}\qquad N=1\\
	\sqrt{M_{2}}\qquad N=2
\end{cases}
 \end{equation}
 and then
 \begin{equation}\label{2.14}
   \left \| u \right \|_{L^{1}B(x,\delta )}\leq (2\delta )^{\frac{N}{2}}M.
 \end{equation}
 Now we proceed to improve the $L^{2}$ boundedness of $u$ to $L^{\infty}$, which is based on the fact that for all $(x,t)\in  \mathbb{R}^{N}\times [0,T_{max})$ and Lemma \ref{lem3}-\ref{lemma2.6}, Proposition \ref{proposition.1}, we have
\begin{equation}
\begin{aligned}
  \nonumber&0\leq u(x,t)\leq \left \| \mathcal{S} _{\alpha }(t)u_{0} \right \|_{L^{\infty }(\mathbb{R}^{N})}+\int_{0}^{t}(t-s)^{\alpha -1}\left \| \mathcal{K}_{\alpha }(t-s)h(s) \right \|_{L^{\infty }(\mathbb{R}^{N})}ds\\
\nonumber&\leq C_{4}\left \| u_{0} \right \|_{L^{\infty }(\mathbb{R}^{N})}+C_{4}\int_{0}^{t}(t-s)^{\alpha -1}\left \| h(s) \right \|_{L^{\infty }(\mathbb{R}^{N})}ds\\
\nonumber&\leq C_{4}\left \| u_{0} \right \|_{L^{\infty }(\mathbb{R}^{N})}+\mu C_{4}\int_{0}^{t}(t-s)^{\alpha -1}\left \| u^{2}(s) \right \|_{L^{\infty }(\mathbb{R}^{N})}ds\\
\nonumber&\leq  C_{4}\left \| u_{0} \right \|_{L^{\infty }(\mathbb{R}^{N})}+\mu C_{4}\int_{0}^{t}(t-s)^{\alpha -1}\left \| u(s) \right \|^{2}_{L^{\infty }(\mathbb{R}^{N})}ds\\
\nonumber&\leq C_{4}\left \| u_{0} \right \|_{L^{\infty }(\mathbb{R}^{N})}+\mu C_{4}M^{2}\int_{0}^{t}(t-s)^{\alpha -1}ds\\
\nonumber&\leq C_{4}\left \| u_{0} \right \|_{L^{\infty }(\mathbb{R}^{N})}+\mu C_{4}M^{2}T^{\alpha }\frac{1}{\alpha}.
\end{aligned}
\end{equation}
So there is
\begin{equation}\label{2.22}
0\leq u(x,t)\leq C_{4}\left \| u_{0} \right \|_{L^{\infty }(\mathbb{R}^{N})}+\mu C_{4}M^{2}T^{\alpha }\frac{1}{\alpha},
\end{equation}
with
\begin{equation}
    M=\begin{cases}
	\sqrt{ [(2\delta \left \| u_{0} \right \|_{L^{\infty }(\mathbb{R}^{N})}^{2})^{1-p}+\frac{(Q+(2\gamma -2C_{2}))T^{\alpha }}{\alpha \Gamma (\alpha )}]^{\frac{1}{1-p}}}, \quad N=1\\
	\sqrt{[((2\delta)^{2}\left \| u_{0} \right \|_{L^{\infty}(\mathbb{R}^{N})}^{2})^{1-p}+\frac{(2\mu C_{GN}^{4}(2,\delta)+(2\gamma -2C_{2}))T^{\alpha }}{\alpha \Gamma (\alpha )}]^{\frac{1}{1-p}}},\quad N=2
\end{cases}
 \end{equation}
 where $Q_{1}=2\mu \left ( (\mu ^{\frac{1}{3}}C_{GN}^{\frac{4}{3}}(1,\delta)+1)^{6}(\eta k)^{-5}+C_{GN}^{10}(1,\delta ) \right )$. From Proposition \ref{Proposition2.1} and inequality \eqref{2.13}, we obtain $T_{max}=+\infty$. The global boundedness of $u$ is obtained in time and the unique classical solution of \eqref{1}-\eqref{2} on $(x,t)\in \mathbb{R}^{N}\times [0,+\infty )$. Theorem \ref{theorem4} is thus proved.
 \end{pot4}

\section{Long time behavior of solutions (the Allee effect)}
To study the long time behavior of solutions for \eqref{1}-\eqref{2}, by \cite{0Global}, we denote $$F(u):=\mu u^{2}(1-kJ*u)-\gamma u.$$ For $0<\gamma<\frac{\mu }{4k}$, there are three constant solutions for $F(u)=0$: $0,a,A$, where
\begin{equation}\label{1.7}
  a=\frac{1-\sqrt{1-4k\frac{\gamma }{\mu }}}{2k},\qquad A=\frac{1+\sqrt{1-4k\frac{\gamma }{\mu }}}{2k},
\end{equation}
 and satisfy $0<\frac{\gamma }{\mu }<a<A$.
 \begin{proposition}\label{proposition3.1}
   Under the assumptions of Theorem \ref{theorem4}, there is \\ $\left \| u(x,t) \right \|_{L^{\infty }(\mathbb{R}^{N}\times (0,+\infty ))}< a$, the function$$H(x,t)=\int _{B(x,\delta )}h(u(y,t))dy$$ with $$h(u)=Aln\left ( 1-\frac{u}{A} \right )-aln\left ( 1-\frac{u}{a} \right )$$ is nonnegative and satisfies
   \begin{equation}\label{3.1}
     \frac{\partial^{\alpha } H(x,t)}{\partial t^{\alpha }}\leq \left | \bigtriangledown u \right |^{p-2}\Delta H(x,t)-D(x,t)
   \end{equation}
   with $$D(x,t)=\frac{1}{2}(A-a)\mu k\int _{B(x,\delta )}u^{2}(y,t)dy.$$
 \end{proposition}

 \begin{pro}
 Denote $k=\left \| u(x,t) \right \|_{L^{\infty }(R^{N}\times [0,\infty ))}$, then $0<k<a$. From the definition of $h(\cdot)$, it is easy to verify that
 \begin{equation}\label{3.2}
   {h}'(u)=\frac{a}{a-u}-\frac{A}{A-u}=\frac{(A-a)u}{(A-u)(a-u)},
 \end{equation}
 and
 \begin{equation}\label{3.3}
  {h}''(u)=\frac{a}{(a-u)^{2}}-\frac{A}{(A-a)^{2}}=\frac{(Aa-u^{2})(A-a)}{(A-u)^{2}(a-u)^{2}}.
 \end{equation}
 Test \eqref{1} by ${h}'(u)\varphi _{\varepsilon }$ with $\varphi _{\varepsilon }(\cdot )\in C_{0}^{\infty }(B(x,\delta )),\varphi _{\varepsilon }(\cdot )\rightarrow 1$ in $B(x,\delta)$ as $\varepsilon \rightarrow 0$. Integrating by parts over $B(x,\delta)$ and 
 \begin{align*}
  \int _{B(x,\delta )}\varphi _{\varepsilon }{h}'(u)\Delta_{p} udy&=\int _{B(x,\delta )}\varphi _{\varepsilon }{h}'(u)div(\left | \bigtriangledown u \right |^{p-2}\bigtriangledown u)dy\\
&=\int _{B(x,\delta )}\varphi _{\varepsilon }{h}'(u)\left | \bigtriangledown u \right |^{p-2}\Delta udy\\
&=\int _{B(x,\delta )}\varphi _{\varepsilon }\left | \bigtriangledown u \right |^{p-2}(\Delta h(u)-{h}''(u)\left | \bigtriangledown u \right |^{2})dy\\
&=\left | \bigtriangledown u \right |^{p-2}\int _{B(x,\delta )}\varphi _{\varepsilon }\Delta h(u)dy-\int _{B(x,\delta )}\varphi _{\varepsilon }{h}''(u)\left | \bigtriangledown u \right |^{p}dy.  
 \end{align*}
Then, we obtain
 $$
 \begin{aligned}
& \int _{B(x,\delta )}{h}'(u)\varphi _{\varepsilon }\partial _{t}^{\alpha }udy\\
&=\int _{B(x,\delta )}\Delta_{p} u{h}'(u)\varphi _{\varepsilon }dy+\int _{B(x,\delta )}[\mu u^{2}(1-kJ*u)-\gamma u]{h}'(u)\varphi _{\varepsilon }dy\\
&= \int _{B(x,\delta )}\left | \bigtriangledown u \right |^{p-2}\Delta h(u)\varphi _{\varepsilon }dy-\int _{B(x,\delta )}{h}''(u)\left | \bigtriangledown u \right |^{p}\varphi _{\varepsilon }dy\\
&\qquad+\int _{B(x,\delta )}[\mu u^{2}(1-kJ*u)-\gamma u]{h}'(u)\varphi _{\varepsilon }dy.
 \end{aligned}
 $$
 From Lemma \ref{lemma2} and Equation \eqref{4}, we get
 $$\int _{B(x,\delta )}{h}'(u)\partial _{t}^{\alpha }udy\geq \int _{B(x,\delta )}\partial _{t}^{\alpha }h(u)dy=\partial _{t}^{\alpha }\int _{B(x,\delta )}h(u)dy.$$
 And by Lemma \ref{lemma2.8}, noticing
 $$
 \begin{aligned}
 \int _{B(x,\delta )}\Delta h(u)dy&=\int _{B(0,\delta )}\Delta h(u(y+x),t)dy \\
 &=\Delta \int _{B(0,\delta )}h(u(x+y),t)dy=\Delta \int _{B(x,\delta )}h(u)dy,
 \end{aligned}
 $$
 taking $\varepsilon \rightarrow 0$, we obtain
 $$
 \begin{aligned}
& \frac{\partial^{\alpha } }{\partial t^{\alpha }} \int _{B(x,\delta )}h(u)dy\\
&\leq \left | \bigtriangledown u \right |^{p-2}\Delta \int _{B(x,\delta )}h(u)dy-\int _{B(x,\delta )}{h}''(u)\left | \bigtriangledown u \right |^{p}dy\\
&\qquad+\int _{B(x,\delta )}[\mu u^{2}(1-kJ*u)-\gamma u]{h}'(u)dy,
 \end{aligned}
 $$
 which is
 \begin{equation}\label{3.4}
 \begin{aligned}
   &\frac{\partial^{\alpha } }{\partial t^{\alpha }}H(x,t)\\
 &\leq \left | \bigtriangledown u \right |^{p-2}\Delta H(x,t)-\int _{B(x,\delta )}{h}''(u)\left | \bigtriangledown u \right |^{p}dy\\
 &\qquad+\int _{B(x,\delta )}[\mu u^{2}(1-kJ*u)-\gamma u]{h}'(u)dy.
\end{aligned}
 \end{equation}
 By \eqref{1.7}, we can get $$\mu u^{2}(1-ku)-\gamma u=k\mu u(A-u)(u-a),$$
and
\begin{equation}\label{3.5}
\begin{aligned}
&\int _{B(x,\delta )}{h}'(u)[\mu u^{2}(1-kJ*u)-\gamma u]dy\\
&=\int _{B(x,\delta )}{h}'(u)[\mu u^{2}(1-ku)-\gamma u]dy+\mu k\int _{B(x,\delta )}{h}'(u)u^{2}(u-J*u)dy\\
&=-(A-a)\mu k\int _{B(x,\delta )}u^{2}(y,t)dy+\mu k\int _{B(x,\delta )}{h}'(u)u^{2}(u-J*u)dy.
\end{aligned}
\end{equation}
Noticing that when $0\leq u\leq k$, there is $$0\leq {h}'(u)u\leq \frac{(A-a)k^{2}}{(A-k)(a-k)}.$$
From Young's inequality and the median value theorem, we can get
\begin{equation}\label{3.6}
\begin{aligned}
&\mu k\int _{B(x,\delta )}{h}'(u)u^{2}(u-J*u)dy\\
&\leq \mu k\int _{B(x,\delta )}\int _{B(x,\delta )}{h}'(u)u^{2}(y,t)(u(y,t)-u(z,t))J(y-z)dzdy\\
&\leq \frac{(A-a)K^{2}}{(A-K)(a-K)}\mu k\int _{B(x,\delta )}\int _{B(x,\delta )}u(y,t)\left | (u(y,t)-u(z,t)) \right |J(y,z)dzdy\\
&\leq \frac{(A-a)K^{4}\mu k}{2(A-K)^{2}(a-K)^{2}}\int _{B(x,\delta )}\int _{B(x,\delta )}(u(z,t)-u(y,t))^{2}J(z-y)dzdy\\
&\qquad+\frac{1}{2}(A-a)\mu k\int _{B(x,\delta )}u^{2}(y,t)dy\\
&\leq \frac{(A-a)K^{4}\mu k}{2(A-K)^{2}(a-K)^{2}}\int _{B(x,\delta )}\int _{B(x,\delta )}\int_{0}^{1}\left | \bigtriangledown u(y+\theta (z-y),t) \right |^{2}\left | z-y \right |^{2}\\
&\qquad J(z-y) d\theta dzdy+\frac{1}{2}(A-a)\mu k\int _{B(x,\delta )}u^{2}(y,t)dy,
  \end{aligned}
\end{equation}
changing the variables ${y}'=y+\theta (z-y) , {z}'=z-y$, then$$\begin{vmatrix}
	& \frac{\partial{y}' }{\partial y} \qquad \frac{\partial {y}'}{\partial z}\\
	& \frac{\partial {z}'}{\partial y}\qquad  \frac{\partial {z}'}{\partial z}\\
\end{vmatrix}=\begin{vmatrix}
	1-\theta   \qquad\theta & \\
	-1\qquad	1 &
\end{vmatrix}=1-\theta+\theta=1.$$
For any $\theta \in  [0,1],y,z\in  B(x,\delta)$, we have ${y}'\in B((1-\theta )x+\theta z,(1-\theta )\delta ) ,{z}'\in B(x-y,\delta )$. Noticing $B((1-\theta )x+\theta z,(1-\theta )\delta )\subseteq B(x,\delta)$ and $B(x-y,\delta)\subseteq B(0,2\delta)$, we obtain
\begin{equation}\label{3.7}
\begin{aligned}
&\frac{(A-a)K^{4}\mu k}{2(A-K)^{2}(a-K)^{2}}\int _{B(x,\delta )}\int _{B(x,\delta )}\int_{0}^{1}\left | \bigtriangledown u(y+\theta (z-y),t) \right |^{2}\left | z-y \right |^{2}\\
& \qquad J(z-y)d\theta dzdy\\
&\leq \frac{(A-a)K^{4}\mu k}{2(A-K)^{2}(a-K)^{2}}\int_{0}^{1}d\theta \int _{B(x,\delta )}\int _{B(x,\delta )}\left | \bigtriangledown u({y}',t) \right |^{2}\left | {z}' \right |^{2}J({z}')d{z}'d{y}'\\
&\leq \frac{(A-a)K^{4}\mu k (2\delta )^{2}}{2(A-K)^{2}(a-K)^{2}}\int _{B(x,\delta )}\left | \bigtriangledown u(y,t) \right |^{2}dy.
\end{aligned}
\end{equation}
Combining \eqref{3.5}-\eqref{3.7}, we obtain
\begin{equation}\label{3.8}
\begin{aligned}
  &\int _{B(x,\delta )}{h}'(u)[\mu u^{2}(1-kJ*u)-\gamma u]dy\\
&\leq -\frac{1}{2}(A-a)\mu k\int _{B(x,\delta )}u^{2}(y,t)dy\\
&\qquad+\frac{(A-a)k^{4}\mu k(2\delta )^{2}}{2(A-a)^{2}(a-k)^{2}}\int _{B(x,\delta )}\left | \bigtriangledown u(y,t) \right |^{2}dy.
\end{aligned}
\end{equation}
From \eqref{3.3}, noticing $0\leq u<a$, we obtain$${h}''(u)> \frac{(A-a)^{2}}{A^{2}a},$$
inserting \eqref{3.8} into \eqref{3.4} and $1<p<2$, we obtain
\begin{equation}\label{3.9}
  \begin{aligned}
  &\frac{\partial^{\alpha } }{\partial t^{\alpha }}H(x,t)\leq \left | \bigtriangledown u \right |^{p-2}\Delta H(x,t)-\frac{(A-a)^{2}}{A^{2}a}\int _{B(x,\delta )}\left | \bigtriangledown u(y,t) \right |^{p}dy\\
  &\qquad+\frac{(A-a)k^{4}\mu k(2\delta )^{2}}{2(A-k)^{2}(a-k)^{2}} \int _{B(x,\delta )}\left | \bigtriangledown u(y,t) \right |^{2}dy-\frac{1}{2}(A-a)\mu k\int _{B(x,\delta )}u^{2}(y,t)dy\\
  &\leq\left | \bigtriangledown u \right |^{p-2} \Delta H(x,t)+[-\frac{(A-a)^{2}}{A^{2}a}+\frac{(A-a)k^{4}\mu k(2\delta )^{2}}{2(A-k)^{2}(a-k)^{2}}]\\
  &\qquad \int _{B(x,\delta )}\left | \bigtriangledown u(y,t) \right |^{2}dy-\frac{1}{2}(A-a)\mu k\int _{B(x,\delta )}u^{2}(y,t)dy.
  \end{aligned}
\end{equation}
By choosing $\delta$ sufficiently small such that$$-\frac{(A-a)^{2}}{A^{2}a}+\frac{(A-a)k^{4}\mu k(2\delta )^{2}}{2(A-k)^{2}(a-k)^{2}}\leq 0,$$ then
$$\frac{\partial^{\alpha } }{\partial t^{\alpha }}H(x,t)\leq \left | \bigtriangledown u \right |^{p-2}\Delta H(x,t)-\frac{1}{2}(A-a)\mu k\int _{B(x,\delta )}u^{2}(y,t)dy,$$
making $$D(x,t)=\frac{1}{2}(A-a)\mu k\int _{B(x,\delta )}u^{2}(y,t)dy.$$
 \end{pro}
 
\begin{remark}
Before proving the global boundedness of equation \eqref{1}-\eqref{2}, we use Laplace's transform to solve the fractional differential equations
\begin{equation}\label{211}
\begin{aligned}
  _{0}^{C}\textrm{D}_{t}^{\alpha }w(t)&=Aw(t)+C_{3},\qquad  0<\alpha <1,t\geq 0\\
 w(0)&=\eta,
 \end{aligned}
\end{equation}
where $_{0}^{C}\textrm{D}_{t}^{\alpha }$ is the Caputo fractional derivative operator \eqref{32}, $A,C_{3}$ is all constants.
\end{remark}
\begin{lemma}\label{lemma2.4}\textsuperscript{\cite{2021zhan}}
  Assume \eqref{211} has a unique continuous solution $w(t),t \in [0,\infty)$, and $A,C_{3}$ is constants, then $\exists T\in [0,\infty )$ such that  $w(t)$ is exponentially bounded in $t \geq T$, which is
  \begin{equation}
    \left \| w\left ( t \right )  \right\| \leq (\left \| \eta \right \|e^{-\sigma T}+\frac{(K+\left | C_{3} \right |)T^{\alpha }e^{-\sigma T}}{\alpha \Gamma (\alpha )})Ce^{[A^{\frac{1}{\alpha }}+\sigma]t}.
  \end{equation}
  \end{lemma}
\begin{lemma}\textsuperscript{\cite{2001Fractional}}
  Assume $E_{\alpha}$ is the Mittag-Leffler type function. Then there is the following inequality
  \begin{equation}\label{2.6}
E_{\alpha }(wt^{\alpha })\leq Ce^{w^{\frac{1}{\alpha }t}},\quad t\geq 0,w\geqslant 0,0<\alpha <2
\end{equation}
where $C$ is a positive constant.
\end{lemma}

\begin{definition}\textsuperscript{\cite{2006Global}}
  Let $X$ be a Banach space,$z_{0}$ belong to $X$, and $f\in L^{1}(0,T;X)$. The function $z(x,t)\in C([0,T];X)$ given by
  \begin{equation}\label{pp}
      z(t)=e^{-t}e^{t\Delta}z_{0}+\int_{0}^{t}e^{-(t-s)}\cdot e^{(t-s)\Delta }f(s)ds,\qquad 0\leq t\leq T,
  \end{equation}
  is the mild solution of \eqref{pp} on $[0,T]$, where $(e^{t\Delta}f)(x,t)=\int _{\mathbb{R}^{N}}G(x-y,t)f(y)dy$ and $G(x,t)$ is the heat kernel by $G(x,t)=\frac{1}{(4\pi t)^{N/2}}exp(-\frac{\left | x \right |^{2}}{4t})$.
\end{definition}
\begin{lemma}\textsuperscript{\cite{2006Global}}\label{lemma5}
  Let $0\leq q\leq p\leq \infty ,\frac{1}{q}-\frac{1}{p}<\frac{1}{N}$ and suppose that $z$ is the function given by \eqref{pp} and $z_{0}\in W^{1,p}(\mathbb{R}^{N})$. If $f\in L^{\infty }(0,\infty ;L^{q}(\mathbb{R}^{N}))$, then
  \begin{eqnarray}
   \nonumber \left \| z(t) \right \|_{L^{p}(\mathbb{R}^{N})}&\leq& \left \|  z_{0}\right \|_{L^{p}(\mathbb{R}^{N})}+C\cdot \Gamma (\gamma )\underset{0<s<t}{sup}\left \| f(s) \right \|_{L^{q}(\mathbb{R}^{N})},\\
\nonumber\left \| \bigtriangledown z(t) \right \|_{L^{p}(\mathbb{R}^{N})}&\leq& \left \|  \bigtriangledown z_{0}\right \|_{L^{p}(\mathbb{R}^{N})}+C\cdot \Gamma (\tilde{\gamma })\underset{0<s<t}{sup}\left \| f(s) \right \|_{L^{q}(\mathbb{R}^{N})},
  \end{eqnarray}
  for $t\in [0,\infty )$, where $C$ is a positive constant independent of $p,\Gamma (\cdot )$ is the gamma function, and $\gamma =1-(\frac{1}{q}-\frac{1}{p})\cdot \frac{N}{2},\tilde{\gamma }=\frac{1}{2}-(\frac{1}{q}-\frac{1}{p})\cdot \frac{N}{2}$.
\end{lemma}

\begin{pot5}
 From the proof of Theorem \ref{theorem4}, for any $\forall{K}'>0$ , from \eqref{2.22} and the definition of
 $M$, there exist $\mu ^{*}({K}')>0$ and $m_{0}({K}')>0$ such that for $\mu \in  (0,\mu ^{*}({K}'))$ and $\left \| u_{0} \right \|_{L^{\infty }(\mathbb{R}^{N})}<m_{0}({K}')$, then for any $\left ( x,t \right )\in \mathbb{R}^{N}\times [0,+\infty ),M$ is sufficiently small such that $$\left \| u(x,t) \right \|_{L^{\infty }(\mathbb{R}^{N}\times [0,+\infty ))}< {K}'.$$
 
 \begin{enumerate}[(i).]
  \item \textbf{The case} $\left \| u(x,t) \right \|_{L^{\infty }(\mathbb{R}^{N}\times [0,+\infty ))}<{K}'=\frac{\gamma }{\mu }$.\\
  \setlength{\parindent}{0pt}
 Noticing $\sigma =\gamma -\mu \left \| u(x,t) \right \|_{L^{\infty }(\mathbb{R}^{N}\times [0,\infty ))}>0$, we have
 \begin{equation}\label{3.10}
   \begin{aligned}
   &\frac{\partial^{\alpha } u}{\partial t^{\alpha }}=\left | \bigtriangledown u \right |^{p-2}\Delta u+\mu u^{2}(1-kJ*u)-\gamma u\\
&=\left | \bigtriangledown u \right |^{p-2}\Delta u+u[\mu u(1-kJ*u)-\gamma ]\\
&=\left | \bigtriangledown u \right |^{p-2}\Delta u-u[\gamma -\mu u(1-kJ*u)]\\
&\leq \left | \bigtriangledown u \right |^{p-2}\Delta u-u[\gamma -\mu u]\\
&\leq\left | \bigtriangledown u \right |^{p-2} \Delta u-u\sigma.
   \end{aligned}
 \end{equation}
 Consider the following fractional ordinary differential equation
 \begin{equation}
   \begin{cases}
   _{0}^{C}\textrm{D}_{t}^{\alpha }w(t)=-\sigma w,\\
	w(0)=\left \| u_{0} \right \|_{L^{\infty }(\mathbb{R}^{N})}.
   \end{cases}
 \end{equation}
 By Lemma \ref{lemma2.4} and Equation \eqref{2.6}, we obtain
 $$u(x,t)\leq w(t)=w(0)E_{\alpha ,1}(-\sigma t)\leq \left \| u_{0} \right \|_{L^{\infty }(\mathbb{R}^{N})}e^{(-\sigma)^{\frac{1}{\alpha }} t}$$
 for all $(x,t)\in \mathbb{R}^{N}\times [0,+\infty )$ and obtain $$\left \| u(\cdot ,t) \right \|_{L^{\infty }(\mathbb{R}^{N})}\leq \left \| u_{0} \right \|_{L^{\infty }(\mathbb{R}^{N})}e^{(-\sigma)^{\frac{1}{\alpha }} t}.$$

    \item \textbf{The case}\, $\left \| u(x,t) \right \|_{L^{\infty }(\mathbb{R}^{N}\times [0,+\infty ))}< {K}'=a$.
 Denoting $H_{0}(y)=H(y,0)$, from \eqref{3.1} in Proposition \ref{proposition3.1} and fractional Duhamel’s formula \ref{lemma2.6}, for any $(x,t)\in \mathbb{R}^{N}\times [0,+\infty )$, we have
$$ H(x,t)\leq \left \| H_{0} \right \|_{L^{\infty }(\mathbb{R}^{N})}-\int_{0}^{t}(t-s)^{\alpha -1}\mathcal{K}_{\alpha }(t-s)D(y,s)ds,$$
from which we obtain$$\int_{0}^{t}(t-s)^{\alpha -1}\mathcal{K}_{\alpha }(t-s)D(y,s)ds \leq \left \| H_{0} \right \|_{L^{\infty }(\mathbb{R}^{N})}.$$
Due to the fact that $u$ is a classical solution, we have that $$\mathcal{K}_{\alpha }(t-s)D(y,s)\in C^{2,1}(\mathbb{R}^{N}\times [0,\infty )),$$
which implies that for all $x \in \mathbb{R}^{N}$, the following limit holds:
$$\lim_{t\rightarrow \infty }\lim_{s\rightarrow t}\mathcal{K}_{\alpha }(t-s)D(y,s)=0,$$
or equivalently
$$\lim_{t\rightarrow \infty }\lim_{s\rightarrow t}\mathcal{K}_{\alpha }(t-s)\int _{B(x,\delta )}u^{2}(z,s)dz=0,$$
which together with the fact that the heat kernel converges to delta function as $s\rightarrow t$, we have that for any $x \in \mathbb{R}^{N}$,
$$\lim_{t\rightarrow \infty }\int _{B(x,\delta )}u^{2}(y,t)dy=0.$$
Furthermore, with the uniform boundedness of $u$ on $\mathbb{R}^{N}\times[0,\infty)$, following Lemma \ref{lemma5}, we can obtain the global boundedness of
$\left \| \bigtriangledown u(\cdot ,t) \right \|_{L^{\infty }(\mathbb{R}^{N})}$, from which and Lemma \ref{lemma2.1} with $\Omega =B(x,\delta ),p=q=\infty ,r=2$, the convergence of $\left \|  u(\cdot ,t)\right \|_{L^{\infty }(B(x,\delta ))}$ follows from the convergences of $\left \|  u(\cdot ,t)\right \|_{L^{2}(B(x,\delta ))}$ immediately. This is, we obtain $$\left \|  u(\cdot ,t)\right \|_{L^{\infty }(B(x,\delta ))}\rightarrow 0$$ as $t \rightarrow 0$.
Therefore, for any compact set in $\mathbb{R}^{N}$, by finite covering, we obtain that $u$ converges to $0$ uniformly in that compact set, which means that $u$ converges locally uniformly to $0$ in $\mathbb{R}^{N}$ as $\rightarrow \infty$. The proof is complete.
\end{enumerate}
\end{pot5} 

 \section{Global bounded of solutions for a nonlinear NTFPLRDE}
\begin{lemma}\label{lemma444}\textsuperscript{\cite{2021zhan}}
   Assume the function $y_{k}(t)$ is nonnegative and exists the Caputo fractional derivative for $t\in [0,T]$ satisfying
   \begin{equation}\label{4.17}
  _{0}^{C}\textrm{D}_{t}^{\alpha }y_{k}(t)\leq -y_{k}+a_{k}(y_{k-1}^{\gamma _{1}}(t)+y_{k-1}^{\gamma _{2}}(t)),
   \end{equation}
   where $a_{k}=\bar{a}3^{rk}>1$ with $\bar{a},r$ are positive bounded constants and $0<\gamma _{2}<\gamma _{1}\leq 3$. Assume also that there exists a bounded constant $K\geq 1$ such that $y_{k}(0)\leq K^{3^{k}}$, then
   \begin{equation}\label{4.18}
    y_{k}(t)\leq (2\bar{a})^{\frac{3^{k}-1}{2}}3^{r (\frac{3^{k+1}}{4}-\frac{k}{2}-\frac{3}{4})}max\left \{ \underset{t\in [0,T]}{sup} y_{0}^{3^{k}}(t),K^{3^{k}}\right \}\frac{T^{\alpha }}{\alpha \Gamma (\alpha)}.
   \end{equation}
 \end{lemma}
 
  \begin{lemma}\textsuperscript{\cite{1959On}}\textsuperscript{\cite{1971Propriet}}\label{lemma33}
   When the parameters $p,q,r$ meet any of the following conditions:

   (i) $q >N \geq 1$, $r\geq1$ and $p=\infty$;

   (ii)$q>max\left \{ 1,\frac{2N}{N+2} \right \},1\leq r<\sigma $ and $r<p<\sigma +1$ in
   \begin{equation}
    \nonumber \sigma :=\begin{cases}
               \frac{(q-1)N+q}{N-q},&q<N,\\
               \infty ,&q\geq N.
              \end{cases}
   \end{equation}
Then the following inequality is established
$$\left \| u \right \|_{L^{p}(\mathbb{R}^{N})}\leq C_{GN}\left \| u \right \|_{L^{r}(\mathbb{R}^{N})}^{1-\lambda ^{*}}
\left \| \bigtriangledown u \right \|_{L^{q}(\mathbb{R}^{N})}^{\lambda ^{*}}$$
among
$$\lambda ^{*}=\frac{qN(p-r)}{p[N(q-r)+qr]}.$$
 \end{lemma}
 
\begin{lemma}\textsuperscript{\cite{Shen2016A}}\label{lemma333}
  Let $N\geq 1$. $p$ is the exponent from the Sobolev embedding theorem, $i.e.$
  \begin{equation}\label{22.1}
    \begin{cases}
      p=\frac{2N}{N-2},& N\geq 3\\
     2<p<\infty,& N=2\\
    p=\infty,& N=1
    \end{cases}
  \end{equation}
  $1\leq r<q<p$ and $\frac{q}{r}< \frac{2}{r}+1-\frac{2}{p}$, then for $v\in {H}'(\mathbb{R}^{N})$ and $v\in L^{r}(\mathbb{R}^{N})$, it holds
  \begin{equation}\label{99.1}
  \begin{aligned}
    \left \| v \right \|_{L^{q}(\mathbb{R}^{N})}^{q}&\leq
     C(N)c^{\frac{\lambda q}{2-\lambda q}}_{0}\left \| v \right \|_{L^{p}(\mathbb{R}^{N})}^{\gamma }+c_{0}\left \| \bigtriangledown u \right \|_{L^{2}(\mathbb{R}^{N})}^{2},\quad N>2,\\
\left \| v \right \|_{L^{q}(\mathbb{R}^{N})}^{q}&\leq C(N)(c^{\frac{\lambda q}{2-\lambda q}}_{0}+c^{-\frac{\lambda q}{2-\lambda q}}_{1})\left \| v \right \|_{L^{r}(\mathbb{R}^{N})}^{\gamma }+c_{0}\left \| \bigtriangledown u \right \|_{L^{2}(\mathbb{R}^{N})}^{2}\\
&\qquad+c_{1}\left \| v \right \|_{L^{2}(\mathbb{R}^{N})}^{2},\quad  N=1,2.
	\end{aligned}
  \end{equation}
  Here $C(N)$ are constants depending on $N,c_{0},c_{1}$ are arbitrary positive constants and
  \begin{equation}\label{112}
  \nonumber \lambda =\frac{\frac{1}{r}-\frac{1}{q}}{\frac{1}{r}-\frac{1}{p}}\in (0,1),\gamma =\frac{2(1-\lambda )q}{2-\lambda q}=\frac{2(1-\frac{q}{p})}{\frac{2-q}{r}-\frac{2}{p}+1}.
  \end{equation}
\end{lemma}

\begin{pot6}
\,
\begin{enumerate}[\bf Step 1.]
\item The $L^{p}$ estimates.
\end{enumerate}
For any $x\in \mathbb{R}^{N}$, multiply \eqref{1.1.4} by $u^{k-1},k>1$ and integrating by parts over $\mathbb{R}^{N}$. From Lemma \ref{lemma1} and Lemma \ref{lemma2}, we have$$u^{k-1}(t)\int _{\mathbb{R}^{N}}(_{0}^{C}\textrm{D} _{t}^{\alpha }u)dx\geq \frac{1}{k}(_{0}^{C}\textrm{D} _{t}^{\alpha }\int _{\mathbb{R}^{N}}u^{k}dx).$$
Thus, we obtain
$$
\begin{aligned}
\frac{1}{k}(_{0}^{C}\textrm{D} _{t}^{\alpha }\int _{\mathbb{R}^{N}}u^{k}dx)&\leq-\left | \bigtriangledown u^{m} \right |^{p-2}\int _{\mathbb{R}^{N}}(u+1)^{m-1}\bigtriangledown u\bigtriangledown (u^{k-1})dx\\
&\qquad+\int _{\mathbb{R}^{N}}u^{k+1}dx(1-\int _{\mathbb{R}^{N}}udx)
-\int _{\mathbb{R}^{N}}u^{k}dx\\
&\leq-\left | \bigtriangledown u^{m} \right |^{p-2}(k-1)\int _{\mathbb{R}^{N}}(u+1)^{m-1}u^{k-2}\left | \bigtriangledown u \right |^{2}dx\\
&\qquad+\int _{\mathbb{R}^{N}}u^{k+1}dx(1-\int _{\mathbb{R}^{N}}udx)-\int _{\mathbb{R}^{N}}u^{k}dx\\
&\leq-\left | \bigtriangledown u^{m} \right |^{p-2}(k-1)\int _{\mathbb{R}^{N}}u^{m+k-3}\left | \bigtriangledown u \right |^{2}dx\\
&\qquad+\int _{\mathbb{R}^{N}}u^{k+1}dx(1-\int _{\mathbb{R}^{N}}udx)
-\int _{\mathbb{R}^{N}}u^{k}dx.
\end{aligned}
$$
Present proof $$\int _{R^{N}}u^{m+k-3}\left | \bigtriangledown u \right |^{2}dx=\frac{4}{(m+k-1)^{2}}\left \| \bigtriangledown u^{\frac{m+k-1}{2}} \right \|_{L^{2}(\mathbb{R}^{N})}^{2}.$$
From
$$
\begin{aligned}
\left \| \bigtriangledown u^{\frac{m+k-1}{2}} \right \|_{L^{2}(\mathbb{R}^{N})}^{2}&=\int _{\mathbb{R}^{N}}(\bigtriangledown u^{\frac{m+k-1}{2}})^{2}dx\\
&=\int _{\mathbb{R}^{N}}\bigtriangledown u^{\frac{m+k-1}{2}}\cdot \bigtriangledown u^{\frac{m+k-1}{2}}dx\\
&=\int _{\mathbb{R}^{N}}\frac{(m+k-1)^{2}}{4}u^{m+k-3}\left | \bigtriangledown u \right |^{2}dx,
\end{aligned}
$$
then the following inequality holds
$$
\begin{aligned}
\frac{1}{k}(_{0}^{C}\textrm{D} _{t}^{\alpha }\int _{\mathbb{R}^{N}}u^{k}dx)&\leq-\left | \bigtriangledown u^{m} \right |^{p-2}\frac{4(k-1)}{(m+k-1)^{2}}\left \| \bigtriangledown u^{\frac{m+k-1}{2}} \right \|_{L^{2}(\mathbb{R}^{N})}^{2}+\int _{\mathbb{R}^{N}}u^{k+1}dx\\
&\quad -\int _{\mathbb{R}^{N}}u^{k+1}dx\int _{R^{N}}udx-\int _{\mathbb{R}^{N}}u^{k}dx.
\end{aligned}
$$
Move the term to simplify and get
\begin{equation}\label{4.1.1}
  \begin{aligned}
  &_{0}^{C}\textrm{D} _{t}^{\alpha }\int _{\mathbb{R}^{N}}u^{k}dx+k\left | \bigtriangledown u^{m} \right |^{p-2}\frac{4(k-1)}{(m+k-1)^{2}}\left \| \bigtriangledown u^{\frac{m+k-1}{2}} \right \|_{L^{2}(\mathbb{R}^{N})}^{2}\\
  &\qquad+k\int _{\mathbb{R}^{N}}u^{k+1}dx\int _{\mathbb{R}^{N}}udx+k\int _{\mathbb{R}^{N}}u^{k}dx\\
&\leq k\int _{\mathbb{R}^{N}}u^{k+1}dx.
  \end{aligned}
\end{equation}
The following estimate $k\int _{\mathbb{R}^{N}}u^{k+1}dx$. when $m\leq3$ and
\begin{equation}\label{5.5}
 k> max\left \{ \frac{1}{4}(3-m)(N-2) -(m-1),(1-\frac{m}{2})N-1,N(2-m)-2\right \},
\end{equation}
from Lemma \ref{lemma33}, we obtain
$$
\begin{aligned}
&k\int _{\mathbb{R}^{N}}u^{k+1}dx=k\left \| u^{\frac{k+m-1}{2}} \right \|_{L^{\frac{2(k+1)}{k+m+1}}(\mathbb{R}^{N})}^{\frac{2(k+1)}{k+m-1}}\\
&\leq k\left \| u^{\frac{k+m-1}{2}} \right \|_{L^{\frac{2(k+1)}{k+m+1}}(\mathbb{R}^{N})}^{\frac{2(k+1)}{k+m-1}\lambda ^{*}}\left \| \bigtriangledown u^{\frac{k+m-1}{2}} \right \|_{L^{2}(\mathbb{R}^{N})}^{\frac{2(k+1)(1-\lambda ^{*})}{k+m-1}}.
\end{aligned}
$$
Knowing from \eqref{5.5}
$$\lambda ^{*}=\frac{1+\frac{(k+m-1)N}{2(k+1)}-\frac{N}{2}}{1+\frac{(k+m-1)N}{k+2}-\frac{N}{2}}\in \left \{ max\left \{ 0,\frac{2-m}{k+1} \right \},1 \right \},$$
using Young's inequality, there are
\begin{equation}\label{4.1.3}
\begin{aligned}
k\int _{\mathbb{R}^{N}}u^{k+1}dx&\leq k\left \| u^{\frac{k+m-1}{2}} \right \|_{L^{\frac{k+2}{k+m+1}}(\mathbb{R}^{N})}^{\frac{2(k+1)}{k+m-1}\lambda ^{*}}\left \| \bigtriangledown u^{\frac{k+m-1}{2}} \right \|_{L^{2}(R^{N})}^{\frac{2(k+1)(1-\lambda ^{*})}{k+m-1}}\\
&\leq \frac{2k(k-1)}{(m+k-1)^{2}}\left \| \bigtriangledown u^{\frac{k+m-1}{2}} \right \|_{L^{2}(R^{N})}^{2}\\
&\qquad+C_{1}(N,k,m)\left \| u^{\frac{k+m-1}{2}} \right \|_{L^{\frac{k+2}{k+m-1}}(\mathbb{R}^{N})}^{Q}
\end{aligned}
\end{equation}
where is $Q=\frac{2(k+1)\lambda ^{*}}{m-2+(k+1)\lambda ^{*}}$.
Next estimate $\left \| u^{\frac{k+m-1}{2}} \right \|_{L^{\frac{k+2}{k+m-1}}(\mathbb{R}^{N})}^{Q}$.
We will use the interpolation inequality to get
\begin{equation}\label{4.1.4}
  \begin{aligned}
  \left \| u^{\frac{k+m-1}{2}} \right \|_{L^{\frac{k+2}{k+m-1}}(\mathbb{R}^{N})}^{Q}&\leq\left \| u^{\frac{k+m-1}{2}} \right \|_{L^{\frac{2(k+1)}{k+m-1}}(\mathbb{R}^{N})}^{Q\lambda }\left \| u^{\frac{k+m-1}{2}} \right \|_{L^{\frac{2}{k+m-1}}(\mathbb{R}^{N})}^{Q(1-\lambda )}\\
&\leq(\left \| u^{\frac{k+m-1}{2}} \right \|_{L^{\frac{2(k+1)}{k+m-1}}(\mathbb{R}^{N})}^{\frac{2(k+1)}{k+m-1} }\left \| u^{\frac{k+m-1}{2}} \right \|_{L^{\frac{2}{k+m-1}}(\mathbb{R}^{N})}^{\frac{2}{k+m-1}})^{\frac{Q\lambda (k+m-1)}{2(k+1)}}\\
&\left \| u^{\frac{k+m-1}{2}} \right \|_{L^{\frac{2}{k+m-1}}(\mathbb{R}^{N})}^{Q(1-\lambda -\frac{\lambda }{k+1})}
  \end{aligned}
\end{equation}
in $\lambda =\frac{k+1}{k+2}$, and $$Q(1-\lambda -\frac{\lambda }{k+1})=0.$$ Then
\begin{equation}\label{4.1.5}
\begin{aligned}
  &C_{1}(N,k,m)\left \| u^{\frac{k+m-1}{2}} \right \|_{L^{\frac{k+2}{k+m-1}}(\mathbb{R}^{N})}^{Q}\\
  &\leq C_{1}(N,k,m)(\left \| u^{\frac{k+m-1}{2}} \right \|_{L^{\frac{2(k+1)}{k+m-1}}(\mathbb{R}^{N})}^{\frac{2(k+1)}{k+m-1} }\left \| u^{\frac{k+m-1}{2}} \right \|_{L^{\frac{2}{k+m-1}}(R^{N})}^{\frac{2}{k+m-1}})^{\frac{Q\lambda (k+m-1)}{2(k+1)}}.
\end{aligned}
\end{equation}
Noticing that when $m>2-\frac{2}{N}$, it is easy to verify$$\frac{Q\lambda (k+m-1)}{2(k+1)}=\frac{(k+1)(m-2)+\lambda ^{*}(k+1)(3-m)}{(k+1)(k+m-1)\lambda ^{*}}< 1,$$
using Young's inequality, then
\begin{equation}\label{4.1.6}
\begin{aligned}
&C_{1}(N,k,m)(\left \| u^{\frac{k+m-1}{2}} \right \|_{L^{\frac{2(k+1)}{k+m-1}}(\mathbb{R}^{N})}^{\frac{2(k+1)}{k+m-1} }\left \| u^{\frac{k+m-1}{2}} \right \|_{L^{\frac{2}{k+m-1}}(\mathbb{R}^{N})}^{\frac{2}{k+m-1}})^{\frac{Q\lambda (k+m-1)}{2(k+1)}}\\
 &\leq k\left \| u^{\frac{k+m-1}{2}} \right \|_{L^{\frac{2(k+1)}{k+m-1}}(\mathbb{R}^{N})}^{\frac{2(k+1)}{k+m-1} }\left \| u^{\frac{k+m-1}{2}} \right \|_{L^{\frac{2}{k+m-1}}(\mathbb{R}^{N})}^{\frac{2}{k+m-1}}+C_{2}(N,k,m).
\end{aligned}
\end{equation}
Substitute \eqref{4.1.3}-\eqref{4.1.6} into \eqref{4.1.1} and $\left | \bigtriangledown u^{m} \right | ^{p-2}<1$, then we will get
\begin{equation}\label{4.1.7}
  _{0}^{C}\textrm{D}_{t}^{\alpha }\int _{\mathbb{R}^{N}}u^{k}dx+k\int _{\mathbb{R}^{N}}u^{k}dx\leq C_{2}(N,k,m),
\end{equation}
from Lemma \ref{corollary1}, we obtain
\begin{equation}\label{4.1.8}
  \int _{\mathbb{R}^{N}}u^{k}dx\leq \left \| u_{0} \right \|_{L^{k}(\mathbb{R}^{N})}^{k}+\frac{C_{2}(N,k,m)T^{\alpha }}{\alpha \Gamma },t\in [0,T].
\end{equation}

\begin{enumerate}[\bf Step 2.]
\item  The $L^{\infty}$ estimates.
\end{enumerate}
On account of the above arguments, our last task is to give the uniform boundedness of solution for any $t>0$.
Denote $q_{k}=2^{k}+2$, by taking $k=q_{k}$ in \eqref{4.1.1}, we have
\begin{equation}\label{4.2.1}
  \begin{aligned}
  &_{0}^{C}\textrm{D}_{t}^{\alpha }\int _{\mathbb{R}^{N}}u^{q_{k}}dx+\left | \bigtriangledown u^{m} \right | ^{p-2}\frac{4q_{k}(q_{k}-1)}{(m+q_{k}-1)^{2}}\left \| \bigtriangledown u^{\frac{m+q_{k}-1}{2}} \right \|_{L^{2}(\mathbb{R}^{N})}^{2}+q_{k}\int _{\mathbb{R}^{N}}u^{q_{k}+1}dx\int _{\mathbb{R}^{N}}udx\\
  &+q_{k}\int _{\mathbb{R}^{N}}u^{q_{k}}dx\leq q_{k}\int _{\mathbb{R}^{N}}u^{q_{k}+1}dx
  \end{aligned}
\end{equation}
armed with Lemma \ref{lemma333}, letting
$$v=\frac{m+q_ {k}-1}{2},q=\frac{2(q_{k}+1)}{m+q_{k}-1},r=\frac{2q_{k-1}}{m+q_{k}-1},c_{0}=c_{1}=\frac{1}{2q_{k}},$$
one has that for $N \geq 1$,
\begin{equation}\label{4.2.2}
\begin{aligned}
  \left \| u \right \|_{L^{q_{k}+1}(\mathbb{R}^{N})}^{q_{k}+1}&\leq C(N)c_{0}^{\frac{1}{\delta _{1}-1}}(\int _{\mathbb{R}^{N}}u^{q_{k}}dx)^{\gamma _{1}}
+\frac{1}{2q^{k}}\left \| \bigtriangledown  u^{\frac{m+q_{k}-1}{2}} \right \|_{L^{2}(\mathbb{R}^{N})}^{2}\\
&\qquad+\frac{1}{2q_{k}}\left \| u \right \|_{L^{m+q_{k}-1}(\mathbb{R}^{N})}^{m+q_{k}-1},
\end{aligned}
\end{equation}
where
$$\gamma _{1}=1+\frac{q_{k}+q_{k-1}+1}{q_{k-1}+\frac{p(m-2)}{p-2}}\leq 2,$$
$$\delta _{1}=\frac{(m+q_{k}-1)-2\frac{q_{k-1}}{p^{*}}}{q_{k}-q_{k-1}+1}=O(1).$$
Substituting \eqref{4.2.2} into \eqref{4.2.1} and with notice that $\frac{4q_{k}(q_{k}-1)}{(m+q_{k}-1)^{2}}\geq 2$ and $\left | \bigtriangledown u^{m} \right | ^{p-2}<1$. It follows
\begin{equation}\label{4.2.3}
  \begin{aligned}
  &_{0}^{C}\textrm{D}_{t}^{\alpha }\int _{\mathbb{R}^{N}}u^{q_{k}}dx+\frac{3}{2}\left \| \bigtriangledown u^{\frac{m+q_{k}-1}{2}} \right \|_{L^{2}(\mathbb{R}^{N})}^{2}+q_{k}\int _{\mathbb{R}^{N}}u^{q_{k}}dx\\
 &\qquad +q_{k}\int _{\mathbb{R}^{N}}u^{q_{k}+1}dx\int _{\mathbb{R}^{N}}udx\\
&\leq q_{k}\int _{\mathbb{R}^{N}}u^{q_{k}+1}dx\leq c_{1}(N)q_{k}^{\frac{\delta _{1}}{\delta _{1}-1}}(\int _{\mathbb{R}^{N}}u^{q_{k-1}}dx)^{\gamma _{1}}+\frac{1}{2}\left \| u \right \|_{L^{m+q_{k}-1}(\mathbb{R}^{N})}^{m+q_{k}-1}.
  \end{aligned}
\end{equation}
Applying Lemma \ref{lemma333} with $$v=u^ {\frac{m+q_{k}-1}{2}},q=2,\gamma=\frac{2q_{k-1}}{m+q_{k}-1},c_{0}=c_{1}=\frac{1}{2}$$ noticing $q_{k-1}=\frac{(q_{k}+1)+1}{2}$, and using Young's inequality, we obtain
\begin{equation}\label{4.2.4}
  \begin{aligned}
  &\frac{1}{2}\left \| u \right \|_{L^{m+q_{k}-1}(\mathbb{R}^{N})}^{m+q_{k}-1}=\frac{1}{2}\int _{\mathbb{R}^{N}}u^{m+q_{k}-1}dx\\
&\leq c_{2}(N)(\int _{\mathbb{R}^{N}}u^{q_{k-1}}dx)^{\gamma _{2}}+\frac{1}{2}\left \| \bigtriangledown u^{\frac{m+q_{k}-1}{2}} \right \|_{L^{2}(\mathbb{R}^{N})}^{2}\\
&\leq \frac{2}{3}\int _{\mathbb{R}^{N}}udx\int _{\mathbb{R}^{N}}u^{q_{k}+1}dx+c_{3}(N)+\frac{1}{2}\left \| \bigtriangledown u^{\frac{m+q_{k}-1}{2}} \right \|_{L^{2}(\mathbb{R}^{N})}^{2},
  \end{aligned}
\end{equation}
where $$\gamma _{2}=1+\frac{m+q_{k}-q_{k-1}-1}{q_{k-1}}< 2.$$
By summing up \eqref{4.2.3} and \eqref{4.2.4}, with the fact that $\gamma _{1}\leq 2$ and $\gamma _{2}<2$, we have
\begin{equation}
\begin{aligned}
&_{0}^{C}\textrm{D}_{t}^{\alpha }\int _{\mathbb{R}^{N}}u^{q_{k}}dx+q_{k}\int _{\mathbb{R}^{N}}u^{q_{k}}dx\leq c_{1}(N)q_{k}^{\frac{\delta _{1}}{\delta _{1}-1}}(\int _{\mathbb{R}^{N}}u^{q_{k-1}}dx)^{\gamma _{1}}+c_{3}(N)\\
&\leq max\left \{ c_{1} (N),c_{3}(N)\right \}q_{k}^{\frac{\delta _{1}}{\delta _{1}-1}}[(\int _{\mathbb{R}^{N}}u^{q_{k-1}}dx)^{\gamma 1}+1]\\
&\leq 2max\left \{ c_{1} (N),c_{3}(N)\right \}q_{k}^{\frac{\delta _{1}}{\delta _{1}-1}}max\left \{ (\int _{\mathbb{R}^{N}}u^{q_{k-1}}dx)^{2},1 \right \}.
\end{aligned}
\end{equation}
 Let $K_{0}=max\left \{ 1,\left \| u_{0} \right \|_{L^{1}(\mathbb{R}^{N})},\left \| u_{0} \right \| _{L^{\infty }(\mathbb{R}^{N})}\right \}$, we have the following inequality for initial data
\begin{equation}\label{4.2.5}
  \int _{\mathbb{R}^{N}}u_{0}^{q_{k}}dx\leq (max\left \{ \left \| u_{0} \right \|_{L^{1}(\mathbb{R}^{N})},\left \| u_{0} \right \| _{L^{\infty }(\mathbb{R}^{N})}\right \})^{q_{k}}\leq K_{0}^{q_{k}}.
\end{equation}
 Let $d_{0}=\frac{\delta _{1}}{\delta _{1}-1}$, it is easy to that $q_{k}^{d_{0}}=(2^{k}+2)^{d_{0}}\leq (2^{k}+2^{k+1})^{d_{0}}$. By taking $\bar{a}=max\left \{ c_{1}(N),c_{3}(N) \right \}3^{d_{0}}$ in the Lemma \ref{lemma444}, we obtain
\begin{equation}\label{4.2.6}
  \int u^{q_{k}}dx\leq (2\bar{a})^{2^{k}-1}2^{d_{0}(2^{k+1}-k-2)}max\left \{ \underset{t\geq 0}{sup}(\int _{\mathbb{R}^{N}}u^{q}dx)^{2^{k}},k_{0}^{q_{k}} \right \}\frac{T^{\alpha }}{\alpha \Gamma (\alpha )}.
\end{equation}
 Since $q_{k}=2^{k}+2$ and taking the power $\frac{1}{q_{k}}$ to both sides of \eqref{4.2.6}, then the boundedness of the solution $u(x,t)$ is obtained by passing to the limit $k\rightarrow\infty$
 \begin{equation}\label{4.2.7}
   \left \| u(x,t) \right \|_{L^{\infty }(\mathbb{R}^{N})}\leq 2\bar{a}2^{2d_{0}}max\left \{ sup_{t\geq 0}\int _{\mathbb{R}^{N}} u^{q_{0}}dx,K_{0}\right \}\frac{T^{\alpha }}{\alpha \Gamma (\alpha )}.
 \end{equation}
 On the other hand, by \eqref{4.1.8} with $q_{0}>2$, we know
 $$\int _{\mathbb{R}^{N}}u^{q_{0}}dx\leq \int _{\mathbb{R}^{N}}u^{3}dx\leq \left \| u_{0} \right \|_{L^{3}(\mathbb{R}^{N})}^{3}+\frac{cT^{\alpha }}{\alpha \Gamma (\alpha )}\leq K_{0}^{3}+\frac{cT^{\alpha }}{\alpha \Gamma (\alpha )}.$$
 Therefore we finally have
 \begin{equation}
   \left \| u(x,t) \right \|_{L^{\infty }(\mathbb{R}^{N})}\leq c(N,\left \| u_{0} \right \|_{L^{1}(\mathbb{R}^{N})},\left \| u_{0} \right \|_{L^{\infty }(\mathbb{R}^{N})},T^{\alpha })=M.
 \end{equation}
 \end{pot6}

\bibliographystyle{elsarticle-num}
\bibliography{reference}

\end{document}